\documentclass[opre]{informs3_revised} 

\OneAndAHalfSpacedXII 

\usepackage{natbib}
 \bibpunct[, ]{(}{)}{,}{a}{}{,}%
 \def\bibfont{\small}%

\usepackage{style_bo}

\TheoremsNumberedThrough     
\ECRepeatTheorems

\EquationsNumberedThrough    

\MANUSCRIPTNO{}
\usepackage{fix-cm}

\def\revision#1{{\color{black}#1}}

\begin{document}



\RUNTITLE{Machine Learning for Bilevel and Stochastic Programming}

\TITLE{Machine Learning-Augmented Optimization of Large Bilevel and Two-stage Stochastic Programs: Application to Cycling Network Design}

\ARTICLEAUTHORS{%
\AUTHOR{Timothy C. Y. Chan, Bo Lin}
\AFF{Department of Mechanical \& Industrial Engineering, University of Toronto, \EMAIL{\{tcychan, blin\}@mie.utoronto.ca}, \URL{}}
\AUTHOR{Shoshanna Saxe}
\AFF{Department of  Civil \& Mineral Engineering, University of Toronto, \EMAIL{s.saxe@utoronto.ca}, \URL{}}
} 

\ABSTRACT{%
A wide range of decision problems can be formulated as bilevel programs with independent followers, which as a special case include two-stage stochastic programs. These problems are notoriously difficult to solve especially when a large number of followers present. Motivated by a real-world cycling infrastructure planning application, we present a general approach to solving such problems.
We propose an optimization model that explicitly considers a sampled subset of followers and exploits a machine learning model to estimate the objective values of unsampled followers. We prove bounds on the optimality gap of the generated leader decision as measured by the original objective function that considers the full follower set. We then develop follower sampling algorithms to tighten the bounds and a representation learning approach to learn follower features, which are used as inputs to the embedded machine learning model. Through numerical studies, we show that our approach generates leader decisions of higher quality compared to baselines. Finally, in collaboration with the City of Toronto, we perform a real-world case study in Toronto where we solve a cycling network design problem with over one million followers. Compared to the current practice, our approach improves Toronto's cycling accessibility by 19.2\%, equivalent to \$18M in potential cost savings. Our approach is being used to inform the cycling infrastructure planning in Toronto and outperforms the current practice by a large margin. It can be generalized to any decision problems that are formulated as bilevel programs with independent followers.
}%



\maketitle

%

\section{Introduction}

This paper is concerned with solving bilevel programs with a large number of followers and where the feasible region of the leader is independent of the followers. \revision{The leader seeks to find a solution that minimizes the total cost of all followers whose decisions depend on both the leader's decision and their own objectives and constraints.} A wide range of decision problems can be modeled this way, including transportation network design \citep{liu2019bike, liu2021urban, lim2021bicycle}, energy pricing \citep{zugno2013bilevel}, and portfolio optimization \citep{carrion2009bilevel, leal2020portfolio}. The main challenge stems from having to model a large number of follower problems to evaluate the leader’s decision. Our method was driven by a real cycling infrastructure planning application in Toronto, Canada that involves over one million followers.

\subsection{Problem Motivation: Cycling Infrastructure Planning} \label{subsec:cycling_app}
Cycling has become an increasingly popular transportation mode due to its positive impact on urban mobility, public health, and the environment. During the COVID-19 pandemic, cycling popularity increased significantly since it represented a low-cost and safe alternative to driving and public transit, while also improving access to essential services \citep{kraus2021provisional}. However, safety and comfort concerns remain major barriers to cycling uptake globally \citep{dill2016revisiting}. Building high-quality cycling infrastructure is among the most effective ways to alleviate cycling stress \citep{buehler2016bikeway}, but its implementation is often constrained by limited financial and human resources. 
In this paper, we develop a bilevel model to optimize the locations of new cycling infrastructure. \revision{This model maximizes a transportation metric called ``low-stress cycling accessibility'', defined as the total amount of ``opportunities'' (e.g., jobs) accessible by individuals via streets that are low-stress (i.e., safe for cycling).} This metric has been shown to be predictive of cycling mode choice \citep{imani2019cycle} and is used to assess existing and new cycling infrastructure in Toronto \citep{city2021plan, city2021active}. In our bilevel model, the leader is a transportation planner who designs a cycling network subject to an infrastructure budget (e.g., 100 km), assuming cyclists will use the low-stress network to travel to opportunities via shortest paths. The followers correspond to all possible origin-destination pairs between units of population and opportunity. The resulting formulation for Toronto includes over one million origin-destination pairs between 3,702 geographic units known as census dissemination areas (DAs). This model is very large, and commercial solvers struggle to find a feasible solution, motivating the development of our method.

\subsection{Technical Challenge} \label{subsec:tech_challenge}

Having a large set of followers $\mS$ adds to the already difficult task of solving bilevel problems as it drastically increases the problem size. 
As we show later, the bilevel problem we consider generalizes two-stage stochastic programming when the leader and follower objectives are identical. As a result, we can draw on approaches from both communities to deal with large $\mS$. Thus, in this paper, readers should think of ``leader'' in a bilevel program and ``first-stage decision maker'' in a two-stage stochastic program as synonymous, and similarly for ``follower'' and ``second-stage decision maker''. As we discuss the bilevel or stochastic programming literature below, we use the corresponding terminology. Two predominant strategies to dealing with large $\mS$ are: (1) solving the problem with a small sample of $\mS$, and (2) approximating follower costs without \revision{explicitly modeling the follower problems}.




Sampling a smaller follower set can be done via random sampling \citep{liu2021urban} or clustering \citep{bertsimas2022optimization}. Given a sample $\mT\subseteq \mS$, a feasible leader solution can be derived \revision{by solving a reduced problem that minimizes the total cost of followers in $\mT$. The reduced problem is easier to solve because its size only depends on $|\mT|$.} However, there is no theoretical guarantee on the performance of the obtained leader solution with respect to the original problem's objective \revision{(i.e., the total cost of followers in $\mS$)}.

For the second strategy, many different algorithms have been developed to approximate the followers' cost. \revision{For example, machine learning (ML) methods have been used to predict the second-stage cost based on the first-stage decision \citep{mivsic2020optimization, liu2021time}. A feasible first-stage solution can then be obtained by solving a surrogate problem that minimizes the output of the trained ML model. Alternatively, the solution of the second-stage problems may be replaced with decision rules that are easier to optimize \citep{chen2008linear}. However, both approaches have limitations. The former method requires the surrogate problem to both be tractable (which often necessitates using simple ML models like linear regression) and achieve strong predictive performance---a combination that is hard to achieve in practice.} On the other hand, using decision rules may lead to infeasible follower solutions when the follower problem includes non-trivial constraints. Moreover, neither method has optimality guarantees.

\revision{In this paper, building on both strategies, we propose an ML-augmented optimization model that is both computationally tractable and capable of generating provably high-quality solutions to the original problem. Specifically, our model explicitly considers a sampled subset $\mT \subseteq \mS$ and augments the objective function with an ML component that approximates the cost of followers in $\mS \backslash \mT$. Compared to pure sampling methods, the ML component captures the broader impact of the leader's decision, resulting in improved solution quality for the original model. Unlike existing approximation methods that rely on pre-trained ML models to map leader decisions to follower costs, our approach embeds the ML model training directly into the bilevel program, so the choice of leader decision influences the trained ML model. This simultaneous optimization and ML model training framework enables the derivation of new theoretical guarantees for the leader's solution.}
\subsection{Contributions} 
\begin{enumerate}
    \item We develop an ML-augmented approach to solving bilevel optimization problems with a large number of independent followers, which generalizes two-stage stochastic programming as a special case (Section \ref{sec:model_pre}). We consider using the $k$-nearest neighbor regression and general parametric regression models (Section \ref{subsec:func_class}). We develop theoretical bounds on the quality of the leader decision from our ML-augmented model as evaluated on the original objective function that considers with the full set of followers (Section \ref{subsec:theo_prop}).
    \item Informed by our theoretical insights, we develop practical strategies to enhance the performance of the ML-augmented model, including i) a follower sampling algorithm to tighten the theoretical bounds (Section \ref{subsec:practical}), and ii) a representation learning method that learns follower features that are predictive of follower objective values (Section \ref{sec:repre_learning}). 
    \item We demonstrate the effectiveness of our approach via computational studies on a cycling network design problem (introduced in Section \ref{subsec:prob_formulation}). We show that i) our learned features are more predictive of follower objective values compared to baseline features from the literature; ii) our follower sampling algorithms further improve the ML models' out-of-sample prediction accuracy by a large margin compared to baseline sampling methods (Section \ref{subsec:exp_pred}); iii) our strong predictive performance translates into high-quality and stable leader decisions from the ML-augmented model. The performance gap between our approach and sampling-based models without the ML component is particularly large when the follower sample is small (Section \ref{subsec:exp_opt}).
    
    \item In collaboration with the City of Toronto, we perform a real-world case study on cycling infrastructure planning (Section \ref{sec:case_study}). We solve a large-scale cycling network design problem where we compare our model against i) purely sampling-based methods that do not use ML and ii) a greedy expansion method that closely matches real-world practice. Compared to i), our method can achieve accessibility improvements between 5.8--34.3\%. Compared to ii), our approach can increase accessibility by 19.2\% on average. If we consider 100 km of cycling infrastructure to be designed using a greedy method, our method can achieve a similar level of accessibility using only 70 km, equivalent to \$18M in potential cost savings. 
\end{enumerate} 

All proofs are in the Electronic Companion.

\section{Literature review} \label{sec:lit_review}


\textbf{Integration of ML and optimization}. There has been tremendous growth in the combination of ML and optimization techniques. ``Predict, then optimize'' is a common modeling paradigm that uses ML models to estimate parameters in an optimization problem, which is then solved with those estimates to obtain decisions \citep{elmachtoub2022smart}. Recent progress has been made in using ML models to prescribe decisions based on contextual features \citep{ban2019big, bertsimas2020predictive}, build end-to-end optimization solvers \citep{khalil2017learningcomb}, and to speed-up optimization algorithms \citep{khalil2016learning, morabit2021machine}. Closest to our work is the literature that integrates pretrained ML models into the solution of optimization problems to map decision variables to uncertain objective values \citep{mivsic2020optimization, liu2021time}. Our work differs in that we integrate the ML model training directly into the optimization problem.

\revision{\textbf{Bilevel optimization.} Bilevel optimization is a modeling framework widely adopted in many areas including transportation \citep{liu2022planning}, energy \citep{carvalho2024nash}, and marketing \citep{li2022convex}. These problems are inherently difficult to solve due to their non-convexity and non-differentiability \citep{colson2007overview}. Classical solution methods typically replace the follower problems with their optimality conditions and then solve the resulting single-level model with algorithms, such as the L-shaped method \citep{birge2011introduction}, the Kuhn-Tucker approach \citep{bard2013practical}, and penalty function methods \citep{white1993penalty}. However, they struggle to handle a large set of followers, as the size of the single-level model scales with the number of followers, which is problematic when the leader's problem is hard, e.g., when it is non-convex. Existing studies have considered reducing the problem size by only modeling a subset of followers, either selected randomly \citep{liu2021urban, lim2021bicycle} or based on domain expertise \citep{zugno2013bilevel}. Despite its promise in improving tractability, the impact of sampling on solution quality remains unexplored. To address this gap, we propose a method that leverages sampling and is complemented by theoretical insights into the relationship between sampling and solution quality. These insights lead to a sampling algorithm that enhances the performance of our approach.}


\textbf{Scenario reduction}. Scenario reduction has been extensively studied in the stochastic programming literature. One stream of literature quantifies the similarity between individual scenarios and then applies clustering methods to select a subset. Common measures include the cost difference between single-scenario optimization problems \citep{keutchayan2021problem}, the opportunity cost of switching between scenarios \citep{bertsimas2022optimization}, and the distance between scenario feature vectors \citep{crainic2014scenario}. 
Another stream selects a scenario subset by minimizing the discrepancy between the distributions described by the two sets as measured using the Wasserstein distance \citep{bertsimas2022optimization} and Fortet-Mourier metrics \citep{dupavcova2003scenario}. 
\revision{More recently, ML has been used to predict the most ``representative'' scenario \citep{bengio2020learning} and construct scenario representations to support scenario clustering \citep{wu2022learning}.} 

\revision{
While our follower selection method is in the spirit of scenario reduction, it has a different goal. Scenario reduction focuses on ``solution stability'', ensuring that the first-stage solution and optimal value given by the reduced model are similar to those of the original model with the full scenario set. In contrast, we emphasize ``solution quality'', ensuring that the first-stage solution from our model is of high quality as evaluated on the original model. This distinction arises because, in our application, the full scenario set represents the actual population of interest, whereas in two-stage stochastic programming, it typically approximates an (unobservable) true distribution. Consequently, solution quality is particularly important in our setting, driving new theoretical analyses and the development of new sampling algorithms.
}

\textbf{Cycling infrastructure planning}. Previous studies on cycling infrastructure planning have considered a variety of approaches. Many papers greedily choose road segments to install cycling infrastructure using expert-defined metrics \citep{olmos2020data}. 
Optimization-based methods typically minimize the travel cost \citep{mauttone2017bicycle}, or maximize total utility \citep{liu2021urban} or ridership \citep{liu2022planning} of a large number of origin-destination (OD) pairs. Due to the large problem size, such models are usually solved heuristically. To the best of our knowledge, only \citet{liu2021urban} and \citet{liu2022planning} solve the problems to optimality at a city scale by randomly sampling OD pairs or restricting the routes that each OD pair can use. Our work adds to the literature by providing a computationally tractable method that can solve larger problems without restrictions such as limited routes for each OD pair. 

\section{Model Preliminaries} \label{sec:model_pre}

In this section, we present the general bilevel problem of interest (Section \ref{subsec:bilevel}), a reduced version based on sampling (Section \ref{subsec:reduced}), and our ML-augmented model (Section \ref{subsec:ml_aug}). 
\revision{We briefly summarize key notation used throughout the paper, which will be restated upon their first introduction. Let $\mS$ denote a set of followers, with $\mT\subseteq \mS$ indicating a sampled subset. The sizes of $\mS$ and $\mT$ are denoted by $m$ and $p$, respectively. For a follower $t\in \mT$, let $m^t$ represent the weight assigned for calculating the ML training loss, and $m^t_{k, \mS \backslash \mT}$ be the number of followers in $\mS \backslash \mT$ whose $k$-nearest neighbors in $\mT$ (measured using a distance metric in the feature space) include follower $t$. We adopt the following notational conventions: vectors and matrices are denoted in bold, sets in calligraphic font, the indicator function by $\mathds{1}(\,\cdot\,)$, $[x]^+ = \max\{0, x\}$, and $[m] = \{1, 2, \dots, m\}$.}

\subsection{The Bilevel Model} \label{subsec:bilevel}
The following is the bilevel optimization problem of interest:
\begin{subequations} \label{prob:general_full}
\begin{align}
    \underset{\bfx, \bfy^1, \dots, \bfy^m}{\textrm{minimize}} \quad 
        & f(\bfx) + \sum_{s\in \mS} q^s g(\bfx, \bfy^s) 
        \label{general:obj_func} \\
    \textrm{subject to} \quad 
        & \bfy^s \in \argmin_{\bfy \in \mY^s(\bfx)} 
            h^s(\bfx, \bfy), \quad \forall s\in \mS\\
        & \bfx \in \mX.
\end{align}
\end{subequations}
Let $\bfx$ denote the leader's decision with a bounded, closed feasible set $\mX \subseteq \bbR^{n_1}$ and cost function $f \colon \bbR^{n_1} \rightarrow \bbR$. Let $\mS$ be a set of $m$ followers, $\bfy^s \in \bbR^{n_2}$ be the decision of follower $s \in \mS$, $g \colon \bbR^{n_1 + n_2} \rightarrow \bbR_+$ measure the cost of a follower's decision, and $q^s \in \bbR_+$ be a nonnegative weight. We assume each follower $s\in \mS$ is optimizing an objective function $h^s\colon \bbR^{n_1 + n_2} \rightarrow \bbR$ subject to a non-empty feasible set $\mY^s(\bfx) \subseteq \bbR^{n_2}$ that depends on the leader's decision. 


\subsubsection{Connection to Two-Stage Stochastic Programs}
    \revision{When $g$ and $h^s$ are identical for all $s\in \mS$, Problem~\eqref{prob:general_full} generalizes a two-stage stochastic program with discrete random variables}
\begin{align*}
    \underset{\bfx, \bfy^1, \dots, \bfy^m}{\textrm{minimize}} \quad 
        & f(\bfx) + \sum_{s\in \mS} q^s g(\bfx, \bfy^s)\\
    \textrm{subject to} \quad 
        & \bfy^s \in \mY^s(\bfx), \quad \forall s\in \mS\\
        & \bfx \in \mX,
\end{align*}
where the decisions of the leader and followers correspond to the first-stage and second-stage decisions, respectively, $\mathcal{S}$ is the set of second-stage scenarios, and $q^s$ (suitably normalized) is the probability of realizing scenario $s$. 

\revision{The benefits of this connection are twofold. First, it inspired the development of Theorems \ref{thm:non-parametric_bound} and \ref{thm:reg_bound}, which characterize the impact of follower sampling on leader solution quality---a type of result commonly seen in the two-stage stochastic programming literature \citep{bertsimas2022optimization, romisch1991stability} but new to bilevel programming. Second, since our bilevel model is more general, the proposed method is also applicable to two-stage stochastic programming, offering a new solution method and theoretical guarantees as discussed in Section \ref{subsec:reduced}.}

\subsubsection{Simplifying the notation.}
To simplify the notation, we write Problem~\eqref{prob:general_full} as
\begin{equation}
    \min_{x\in\mX}F(\bfx),
    \label{prob:general}
\end{equation}
where
\begin{equation}
    F(\bfx) := f(\bfx) + \sum_{s\in\mS}q^sG^s(\bfx)
\end{equation}
and
\begin{equation}
    G^s(\bfx) := \min_{\bfy^s} 
        \left\{ g(\bfx, \bfy^s) \,\middle|\, 
        \bfy^s \in \argmin_{\bfy \in \mY^s(\bfx)} \revision{h^s(\bfx, \bfy)}  \right\},
        \quad \forall s\in \mS.
    \label{eq:follower_cost}
\end{equation}
Let $\bfx^* \in \argmin_{x\in\mX}F(\bfx)$ be an optimal solution to Problem~\eqref{prob:general}. 
We assume $G^s(\bfx)$ is bounded for any $s\in \mS$ and $\bfx \in \mX$. This is a standard condition that renders Problem~\eqref{prob:general} well-defined. Without loss of generality, let $\bar{G} \in \bbR_+$ be a constant such that $0 \leq G^s(\bfx) \leq \bar{G}$ for any $s\in \mS$ and $x\in \mX$.

\subsection{Reduced Model} \label{subsec:reduced}
A common idea to improve the computational tractability of Problem~\eqref{prob:general} is to sample a subset of $\mS$. Given a sampled follower set $\mT\subseteq \mS$, we consider the following \textit{reduced problem}
\begin{equation}
    \min_{x\in\mX} 
    \bar{F}_\mT(\bfx),
    \label{prob:reduced}
\end{equation} 
where
\begin{equation}
   \bar{F}_\mT(\bfx) :=  f(\bfx) + \sum_{t\in \mT} r^t G^t(\bfx),
    \label{prob:naive_sample}
\end{equation}

$G^t(\bfx)$ is as defined in \eqref{eq:follower_cost}, and the weight assigned to scenario $t\in \mT$, $r^t \in \bbR_+$, may be different from $q^t$, due to re-weighting. Let $\bar{\bfx}_\mT$ be an optimal solution to Problem~\eqref{prob:reduced}.  

For two-stage stochastic programming, stability results have been established for Problem~\eqref{prob:reduced}. For example, it is possible to bound $|F(\bfx^*) - \bar{F}_{\mT}(\bar{\bfx}_{\mT})|$  \citep{bertsimas2022optimization} and $\|\bfx^* - \bar{\bfx}_\mT \|$ \citep{romisch1991stability}. However, bounds on $|F(\bfx^*) - F(\bar{\bfx}_{\mT})|$, which we develop in this paper, have only recently been studied \citep{zhang2023optimized} and in a more restricted setting.

\subsection{ML-Augmented Model} \label{subsec:ml_aug}
Given a sampled follower set $\mT \subseteq \mS$, we propose the following \textit{ML-augmented model}
\begin{subequations} \label{prob:pwo_ext}
\begin{align}
    \underset{\bfx \in \mX, P \in \mP}{\textrm{minimize}} \quad
        & f(\bfx) 
        + \sum_{t\in \mT} q^t G^t(\bfx)
        + \sum_{s\in \mS\backslash\mT} q^s P\left(\bff^s\right) 
        \label{pwo:obj_func} \\
    \textrm{subject to} \quad
        \label{pwo:training_loss}
        & \sum_{t\in\mT} 
            m^t\left|P(\bff^t) - G^t(\bfx) \right| \leq \bar{L}.
\end{align}
\end{subequations}
We augment the reduced model by integrating an ML model $P\colon \bbR^\xi \rightarrow \bbR$ that predicts the cost of follower $s\in \mS$ based on a feature vector $\bff^s\in \bbR^\xi$. We use $\mP$ to denote the function class of ML models, and $\bar{L} \in \bbR_+$ to indicate a user-defined upper bound on the training loss. 
\revision{A user-defined weight $m^t$ is assigned to each follower $t \in \mT$ when calculating the training loss.  
}
The training of $P$ on dataset $\left\{\bff^t, G^t(\bfx) \right\}_{t\in \mT}$ is embedded into the problem via the training loss constraint \eqref{pwo:training_loss}. When the ML model can be compactly represented, Problem~\eqref{prob:pwo_ext} only adds a small number of decision variables and constraints to the reduced model (detailed in Section \ref{sec:opt_pred}). Nevertheless, it is expected to generate better leader decisions than the reduced model because the ML prediction helps to capture the impact of leader decisions on out-of-sample followers in $\mS\backslash \mT$. \revision{Additionally, Problem~\eqref{prob:pwo_ext} typically retains the structure of the original bilevel model. Thus, algorithms designed for the original model may still be applicable, further improving the tractability of Problem~\eqref{prob:pwo_ext}.}

\begin{remark}[Choice of Loss function]
    \revision{In Problem~\eqref{prob:pwo_ext}, we use the $L_1$ loss because the resulting model can be easily linearized. However, using other loss functions is also possible and would not invalidate our theoretical results. For example, with the $l_2$ loss, similar solution quality guarantees (Theorem \ref{thm:reg_bound}) can be derived by exploiting the fact that the $l_1$ loss can be upper-bounded by a function of the $l_2$ loss, which follows directly from the Cauchy-Swartz inequality. 
    }     
\end{remark}

\begin{remark}[Relationship between $P$ and $\bfx$]
    During the solution of Problem~\eqref{prob:pwo_ext}, the ML model $P$ can be adjusted on-the-fly according to the leader decision $\bfx$ since both $P$ and $\bfx$ are decision variables. These decisions are connected through constraint~\eqref{pwo:training_loss}, which requires the ML model to achieve a low training loss for any leader decisions. It is necessary to embed the training of $P$ into Problem~\eqref{prob:pwo_ext} because, otherwise, the ML-model prediction would become a constant and thus Problem~\eqref{prob:pwo_ext} would be equivalent to the reduced model. We emphasize that $P$ is not pre-trained. It is trained simultaneously as $\bfx$ is optimized by optimizing the weights in model $P$. 
\end{remark}

\begin{remark}[Relationship between $\bff$ and $\bfx$]
\revision{
    Problem~\eqref{prob:pwo_ext} can be generalized by including $\bfx$ as an input to $P$, replacing $P(\bff)$ with $P(\bff, \bfx)$. When $P$ is non-parametric (e.g., $k$NN), the resulting formulation is equivalent to the current model, as the $\bfx$ component is identical across all followers, and the distance between followers is determined solely by the $\bff$ component. When $P$ is parametric (e.g., linear regression), incorporating $\bfx$ enhances the sensitivity of cost predictions to changes in $\bfx$, enabling it to better distinguish subtle differences between leader decisions. In contrast, in our current model, the ML prediction may be the same for similar leader decisions as the ML model trained for one leader decision may satisfy the training loss constraint under the other. However, the primary drawback of using $P(\bff, \bfx)$ is the introduction of bilinear terms, i.e., the product of ML model weights and leader decisions, which would negatively affect the tractability of Problem~\eqref{prob:pwo_ext}. Given the emphasis on computational efficiency in our real-world application, we focus on ML models that only take $\bff$ as input, leaving the incorporation of $\bfx$ for future research. Despite this simplification, our models still demonstrate strong performance compared to baselines (Sections \ref{subsec:exp_opt} and \ref{sec:case_study}).
}
\end{remark}


Problem~\eqref{prob:pwo_ext} provides a general structure for our modeling approach. Its effectiveness depends on multiple factors: i) function class $\mP$, ii) weighting scheme $m^t$ and upper bound $\bar{L}$, iii) sample $\mT$, and iv) availability of predictive follower features $\bff^s, s \in \mS$. We address the first three items in Section \ref{sec:opt_pred} and the fourth in Section \ref{sec:repre_learning}.

\section{Integrating a Prediction Model} \label{sec:opt_pred}
In Section \ref{subsec:func_class}, we introduce two classes of prediction models -- one non-parametric ($k$NN) and one parametric -- that are compatible with our ML-augmented model. We provide theoretical bounds on performance in Section \ref{subsec:theo_prop}. Finally, we present algorithms and discuss practical implementation, based on insights from examining the bounds, in Section \ref{subsec:practical}.

\subsection{Function Classes} \label{subsec:func_class}

\subsubsection{\texorpdfstring{$k$}{Lg}-nearest neighbor regression.} 
For any fixed $\bfx$, let
\begin{equation}\label{mdl:non-parametric}
    P(\bff^s) = 
        \frac{1}{k}
        \sum_{t\in \mT_k(\bff^s)}
        G^t(\bfx)
\end{equation}
where $k$ denotes the the neighborhood size and $\mT_k(\bff^s) \subseteq \mT$ contains the $k$-nearest neighbors of follower $s$ in the sampled set $\mT$. Since the $k$NN regression does not require training, we do not need constraint~\eqref{pwo:training_loss}. Equivalently, we can simply set $\bar{L} = \infty$. Then, Problem~\eqref{prob:pwo_ext} becomes
\begin{equation}
    \underset{\bfx \in \mX}{\textrm{minimize}} \quad
        f(\bfx) 
        + \sum_{t\in \mT} q^{t} G^{t}(\bfx)
        + \sum_{s\in \mS\backslash\mT}
            \sum_{t\in \mT_k(\bff^s)} 
                \frac{q^{s}}{k}
                G^t(\bfx).
\label{prob:non-parametric_approx}
\end{equation}
Note that our ML-augmented model is also compatible with other non-parametric prediction models, such as locally weighted regression \citep{cleveland1988locally} and kernel regression \citep{parzen1962estimation}, which assign non-uniform weights to nearest neighbors. However, these models would require training, which necessitates including the ``weights'' as decision variables in Problem~\eqref{prob:pwo_ext}, leading to bilinear terms (i.e., the product of ``weights'' and $G^t(\bfx)$) that hinder the model's computational tractability. Since our main focus is to improve the bilevel model's computational tractability, we leave the integration of more sophisticated non-parametric models for future work.



\subsubsection{Parametric regression.} \label{subsubsec:ml_para_reg} Consider a parametric regression model $P(\bff^s; \bstheta)$ parameterized by $\bstheta \in \bsTheta$. Then, Problem~\eqref{prob:pwo_ext} becomes
\begin{equation} \label{prob:reg_approx}
    \min_{\bfx\in \mX, \bstheta \in \bsTheta} 
    \left\{
        f(\bfx) 
        + \sum_{t \in \mT} q^t G^t(\bfx) 
        + \sum_{s\in \mS \backslash \mT}  q^s P(\bff^s; \bstheta) 
    \,\middle|\,
        \sum_{t\in \mT} 
            m^t_{1, \mS\backslash\mT}
            \left| G^t(\bfx) - P(\bff^t; \bstheta) \right| 
        \leq \bar{L}
    \right\},
\end{equation}
where $m^t_{1, \mS\backslash\mT}$ is the number of followers in $\mS\backslash\mT$ whose nearest neighbor in $\mT$ is $t$. 
\revision{
While users may choose any weighting scheme for calculating the training loss, in this paper, we set $m^t = m^t_{1, \mS\backslash\mT}$ to ensure that our theoretical results in Section  \ref{subsec:theo_prop} hold. The intuition is that a training data point should receive a higher weight if it has more ``similar'' test data nearby. This principle has been applied in many ML contexts, e.g., domain adaptation \citep{kouw2019review}. 
}

For Problem~\eqref{prob:reg_approx} to be effective, one should choose a function class that can be compactly represented with $\bstheta$ and $\bff$. For example, a linear regression model $P(\bff; \bstheta) = \bstheta^\intercal \bff$ can be incorporated using only $\xi$ additional continuous decision variables $\bstheta \in \bbR^\xi$. 
An additional set of $|\mT|$ variables and $2|\mT| + 1$ linear constraints are needed to linearize the $L_1$ training loss. Such a representation is roughly in the same complexity class as the reduced model~\eqref{prob:reduced} when $\xi$ and $|\mT|$ are small.
\subsection{Theoretical Properties}
\label{subsec:theo_prop}

\subsubsection{Prediction model setup.} 

We start by formally defining the prediction problem embedded in our ML-augmented model. For any \textit{fixed} leader decision $\bfx$, we are interested in predicting 
a follower's cost $G^s(\bfx)$ based on its features $\bff^s$. We define this regression problem in a feature space $\mF \subseteq \bbR^\xi$ and a target space $\mG_{\bfx} \subseteq \bbR$. We denote by $\eta_{\bfx}\left(\cdot \,|\, \bff \right)$ the probability density function of the target variable given a feature vector $\bff$. We regard $G^s(\bfx)$ as a random variable because the true mapping from features to this target may not be deterministic. For example, consider a network design problem where the follower's cost is the length of the shortest path from an origin to a destination using the network designed by the leader. If we use a one-dimensional binary feature that is 1 if both the origin and destination are in downtown and 0 otherwise, then all downtown OD pairs share the same feature value but with drastically different shortest path lengths. 

\subsubsection{Assumptions.} \label{subsubsec:asm}
Next, we introduce several assumptions that enable the derivation of our theoretical results in Sections \ref{subsubsec:bound_knn} and \ref{subsubsec:bound_reg}.

\begin{assumption}
\label{asm:iid}
For any followers $s, s' \in \mS, s \ne s'$ and leader decisions $\bfx_1, \bfx_2 \in \mX$, the target (random) variables with distributions $\eta_{\bfx_1}(\cdot\,|\, \bff^s) $ and $\eta_{\bfx_2}(\cdot \,|\, \bff^{s'})$ are independent.
\end{assumption}




\begin{assumption}
There exists a constant $\mu \in \bbR_+$ such that,
for any fixed leader decision $\bfx \in \mX$, $\bbE_{G \sim \eta_\bfx(\cdot \,|\, \bff)} \left[ G \,|\, \bff \right]$ is $\mu$-Lipschitz continuous with respect to $\bff$.
\label{asm:lipschitz_cost}
\end{assumption}

\begin{assumption}
\label{asm:lipschitz_pred}
There exists a constant $\lambda \in \bbR_+$ such that, for any fixed ML model parameters $\bstheta \in \bsTheta$, $P(\bff; \bstheta)$ is $\lambda$-Lipschitz continuous with respect to $\bff$.
\end{assumption}

Assumption \ref{asm:iid} implies that  the follower set $\mS$ is independently sampled. This assumption is standard in the ML literature and holds for various applications.
For example, in transportation network design, OD pairs (followers) are usually independently sampled from survey data or ridership data \citep{liu2021urban, liu2022planning}. In two-stage stochastic programming, second-stage scenarios are usually taken from historical observations that can be regarded as independent samples \citep{birge2011introduction}. 
Assumption \ref{asm:lipschitz_cost} limits the change in the expected follower cost as a function of the change in feature space. 
In other words, similar followers (as measured by the distance between their feature vectors) should have similar costs under any leader decisions.
Similar assumptions are commonly made to derive stability results for two-stage stochastic programming where the realized uncertain parameters are used as follower features; See, for example, Assumption 4 in \citet{bertsimas2022optimization} and Theorem 1 in \citet{zhang2023optimized}. Assumption \ref{asm:lipschitz_pred} limits the complexity of $P$, which is critical to avoid overfitting
since $P$ is trained on a small dataset ($\mT$). This condition can be enforced by adding regularization constraints to $\bsTheta$. For example, for linear regression, we can set $\bsTheta = \{ \bstheta \in \bbR^{\xi} \,|\, \|\bstheta\|_1 \leq \lambda \}$. This assumption is needed only for parametric regression models.

Next, we present theoretical bounds for the quality of leader decisions from the $k$NN (Section~\ref{subsubsec:bound_knn}) and parametric regression-augmented models (Section~\ref{subsubsec:bound_reg}), and then follower selection methods that tighten the bounds (Section~\ref{subsec:practical}).

\subsubsection{Bound on \texorpdfstring{$k$}{Lg}NN-augmented model solution.} \label{subsubsec:bound_knn}

\begin{theorem}[Bound on $k$NN-augmented model solution] \label{thm:non-parametric_bound}
Given a follower sample $\mT\subseteq \mS$ and a neighborhood size $k\in \{1, 2, \ldots, |\mT|\}$, let $\bfx^{k\textrm{NN}}_{\mT}$ be an optimal solution to Problem \eqref{prob:non-parametric_approx}, $d_{\mF}$ be a distance metric in $\mF$, $\bar{Q} = \max_{s \in \mS\backslash\mT}q^s$, and $m^t_{k, \mS\backslash\mT} = \sum_{s\in \mS\backslash\mT} \mathds{1}\left[t\in \mT_k(\bff^s) \right]$. If Assumptions \ref{asm:iid} and \ref{asm:lipschitz_cost} hold, then, with probability at least $1 - \gamma$, $F(\bfx^{k\textrm{NN}}_{\mT}) - F(\bfx^*) \leq E_m^{k\textrm{NN}}(\mT)$
where 
\begin{equation*}
    E_m^{k\textrm{NN}}(\mT)
    = 
    \sum_{s\in \mS \backslash \mT}
    \sum_{t\in \mT_k(\bff^s)}
        \frac{2 \mu \bar{Q}}{k}
        d_{\mF}(\bff^s, \bff^t)
    + \sqrt{
        4\bar{Q}^2 \bar{G}^2 
        \left[
            |\mS\backslash\mT| 
            + \sum_{t\in\mT}
                \left(\frac{m^t_{k, \mS\backslash\mT}}{k}\right)^2
        \right]
        \log(1/\gamma)}.
\end{equation*}
\end{theorem}
Theorem \ref{thm:non-parametric_bound} bounds the optimality gap of the solution from the $k$NN-augmented model on the original problem. The first term corresponds to the prediction bias and the second term corresponds to the variance. The first term is proportional to the sum of the distances from each $\bff^s$ to its $k$ nearest neighbors in $\mT$. When the sample size $|\mT|$ is fixed, the second term is controlled by $m^t_{k, \mS\backslash\mT}$. Note that $\sum_{t\in\mT} m^t_{k, \mS\backslash\mT} = |\mS\backslash\mT|$, so the second term is minimized when the $m^t_{k, \mS\backslash\mT}, t \in \mT$ are identical, which follows from the Cauchy-Schwarz inequality. The intuition is that if the followers in $\mS\backslash\mT$ are evenly assigned to sample followers in $\mT$, then the overall prediction performance on $\mS\backslash\mT$ is less affected by the random deviation of the individual cost of follower $t$, $G^t(\bfx)$, from its expected value. We note that $E_m^{k\textrm{NN}}(\mT) = 0$ when $\mT = \mS$ because in this case $\mS \backslash \mT = \emptyset$ and $m^t_{k, \mS\backslash\mT} = 0$.


\subsubsection{Bound on parametric regression-augmented model solution.}
\label{subsubsec:bound_reg}

\begin{theorem} \label{thm:reg_bound}
Given a follower sample $\mT \subseteq \mS$,  $\bfx^{\textrm{PR}}_{\mT}$ be the optimal solution to Problem~\eqref{prob:reg_approx}, $\nu(s)$ be the nearest neighbor of $\bff^s$ in $\{\bff^t\}_{t\in \mT}$, and $m^t_{1, \mS\backslash \mT} = \sum_{s\in \mS\backslash \mT} \mathds{1}[\nu(s) = t]$. If Assumptions \ref{asm:iid}--\ref{asm:lipschitz_pred} hold, with probability at least $1 - \gamma$, $F(\bfx^\text{PR}_{\mT}) - F(\bfx^*) \leq E_m^\textrm{PR}(\mT, \bar{L})$ where
\begin{equation*}
    E_m^\textrm{PR}(\mT, \bar{L})
    = 
    2 \bar{Q} \bar{L} 
    \mathds{1}(\mT \subset \mS)
    + 2 \bar{Q} (\lambda + \mu)
    \sum_{s\in \mS\backslash\mT} 
        d_{\mF}(\bff^s, \bff^{\nu(s)})
    + \sqrt{4\bar{Q}^2\bar{G}^2\left[|\mS\backslash\mT| + \sum_{t\in\mT}(m^t_{1, \mS\backslash\mT})^2 \right]\log(1/\gamma)}.
\end{equation*}
\end{theorem}

Theorem \ref{thm:reg_bound} bounds the optimality gap of the leader's solution from Problem~\eqref{prob:reg_approx} on the original problem. The first term is controlled by the training loss $\bar{L}$, while others are controlled by $\mT$. To reduce the last two terms, $\mT$ should be chosen such that followers $s \in \mS\backslash\mT$ are not too far from its nearest neighbor in $\mT$ (second term) and the assignment of followers in $\mS\backslash\mT$ to followers in $\mT$ should be even (third term). Similar to $E^{k\textrm{NN}}_{m} (\mT)$, $E^{\textrm{PR}}_{m} (\mT)$ converges to zero as $\mT$ goes to $\mS$.

\subsection{Practical Implementation} \label{subsec:practical}

\subsubsection{\texorpdfstring{$k$}{Lg}NN-augmented model.} \label{subsubsec:prac_non-parametric}

Theorem \ref{thm:non-parametric_bound} characterizes the impact of follower selection on the quality of the leader's decision from Problem~\eqref{prob:non-parametric_approx}. While one might be tempted to select $\mT$ by directly minimizing $E_m^{k\textrm{NN}}(\mT)$, solving this problem is challenging due to the complex function form of the variance term. Instead, we propose to select $\mT$ by minimizing the more tractable bias term with constraints that aid in reducing the variance term. Finally, we justify our sample selection by demonstrating the tightness of our bound from Theorem~\ref{thm:non-parametric_bound} when our follower sample is used.

Specifically, we select follower samples by solving the following problem 
\begin{equation}
    \mT^k_{p, d} := 
    \underset{\mT\subseteq \mS}{\argmin}
    \left\{
        \sum_{s\in \mS \backslash \mT}
        \sum_{t\in \mT_k(\bff^s)}
        d_{\mF}(\bff^s, \bff^t)
    \,\middle|\,
        |\mT| \leq p,
        |\mS(t)| \leq d, \forall t\in \mT
    \right\}
    \label{prob:follower_select_knn},
\end{equation}
where $\mS(t) = \{s\in \mS\backslash\mT\, |\, t\in \mT_{k}(\bff^s)\}$ denotes the set of unsampled followers that are assigned to the sampled follower $t\in \mT$ according to the $k$-nearest neighbor rule, $p$ is an upper bound on the sample size imposed by the available computational resources, and $d \geq \lceil m/p \rceil$ is a finite positive constant. The goal of solving Problem~\eqref{prob:follower_select_knn} is to minimize the bias term in $E^{k\textrm{NN}}_{m}(\mT)$. The constraints $|\mS(t)| \leq d$ for $t\in \mT$ induce more even assignment of the unsampled followers to the sampled followers so that the variance term is reduced. 

Next, we demonstrate the tightness of our bound in Theorem \ref{thm:non-parametric_bound} when $\mT^k_{p, d}$ is used. We focus on analyzing the case of $k=1$ because the sample size $p$ is usually small due to the bilevel structure of the $k$NN-augmented model, implying that the optimal choice of $k$ is likely small \citep{stone1977consistent, bickel1983sums}. This conjecture is validated by our empirical analysis based on cross-validation, which consistently yields the choice of $k=1$.

\begin{theorem} \label{thm:knn_tightness}
    If $(\bff^1, \bff^2, \dots, \bff^m)$ is a sequence of i.i.d random vectors in $[0, 1]^\xi$ following a continuous density function, $\xi\geq 2$, $p = \max\{1, \alpha m^{(\xi-1)/\xi}$\} for some $\alpha \in (0, 1]$, $\lceil m/p \rceil \leq d \leq \beta \lceil m/p \rceil$ for some $\beta \geq 1$, and $\gamma \in (0, 1]$, then    \begin{equation*}
        \lim_{m \rightarrow \infty} 
            \frac{1}{m}
            E_m^{1\textrm{NN}}(\mT^1_{p, d})
            = 0.
    \end{equation*}
\end{theorem}

Theorem \ref{thm:knn_tightness} states that even if $E^{k\textrm{NN}}_m(\mT)$ involves the sum of $|\mS\backslash \mT|$ terms, it grows sub-linearly in $|\mS|$, which sheds light on the tightness of the bound. While Theorem \ref{thm:non-parametric_bound} holds for any $\mT \subseteq \mS$, Theorem \ref{thm:knn_tightness} holds only for $\mT^1_{p,d}$, highlighting the importance of intelligent sample selection. Theorem~\ref{thm:knn_tightness} requires the sample size to increase at a rate of $m^{(\xi - 1)/\xi}$, which hints at the practical need to consider larger samples for larger problems. We show in Section \ref{sec:comp} that increasing the sample size generally leads to decisions of higher quality. The rate of increase depends on $\xi$. We thus should use compact follower features whenever possible. If high-dimensional features ($\xi$ large) are necessary, the sample size would grow almost linearly in $|\mS|$. But we would still expect a significant improvement in computation time compared to solving the original problem when the leader/follower problems are non-convex, since computation time would grow exponentially in problem size. Finally, the assumption of $\bff^1, \bff^2\, \ldots, \bff^m$ being in $[0, 1]^\xi$ is nonrestrictive as we can create this structure by applying the min-max standardization to follower features.

\subsubsection{Parametric regression.} \label{subsubsec:prac_reg}

The bound in Theorem \ref{thm:reg_bound} is controlled by i) the sample $\mT$ and ii) the value of $\bar{L}$. To select $\mT$, we consider minimizing the bias term in $E^\textrm{PR}_m(\mT, \bar{L})$ for reasons discussed in the previous section. Since the bias terms in $E^{\textrm{PR}}_m(\mT, \bar{L})$ and $E^{1\textrm{NN}}_m(\mT)$ are identical, we also use $\mT^1_{p,d}$ as defined in Problem~\eqref{prob:follower_select_knn} as our follower sample. For $\bar{L}$, choosing a small value will reduce the bound, but could lead to overfitting or even worse, render Problem~\eqref{prob:reg_approx} infeasible. We view $\bar{L}$ as a hyperparameter that should be tuned and provide an approach for doing so in \ref{subapp:l_bar}.

\begin{theorem} \label{thm:reg_tightness}
    If $(\bff^1, \bff^2, \dots, \bff^m)$ is a sequence of i.i.d. random vectors in $[0, 1]^\xi$ following a continuous density function, where $\xi \geq 2$, $p = \max\{1, \alpha m^{(\xi-1)/ \xi}$\} for some $\alpha \in (0, 1]$, and $\lceil m/p \rceil \leq d \leq \beta \lceil m/p \rceil$ for some $\beta \geq 1$, $\bar{L}$ is a finite positive constant, then
    \begin{equation*}
        \lim_{m \rightarrow \infty} 
            \frac{1}{m}
            E_m^{\textrm{PR}}(\mT^1_{p, d}, \bar{L})
            = 0.
    \end{equation*}
\end{theorem} 

Similar to Theorem \ref{thm:knn_tightness}, Theorem \ref{thm:reg_tightness} comments on the tightness of the bound in Theorem~\ref{thm:reg_bound} when the follower samples are selected using our approach. 

\section{Computational Study: Algorithm Performance on Synthetic Cycling Network Design Problem} \label{sec:comp}

In this section, we validate the effectiveness of our ML-augmented model with our representation learning framework on a cycling network design problem. Computational results in this section are generated using synthetic problem instances for which the original bilevel model~\eqref{prob:general} can be solved to optimality using a Benders decomposition method, allowing accurate evaluation of the reduced model~\eqref{prob:reduced} and our ML-augmented models. We introduce the problem and its formulation in Section \ref{subsec:prob_formulation}, \revision{followed by a method for learning follower features in Section \ref{sec:repre_learning}}. We present two experiments to validate the predictive power of the learned follower features and the value of integrating an ML model in the optimization problem in Sections \ref{subsec:exp_pred} and \ref{subsec:exp_opt}, respectively. 

\subsection{Maximum Accessibility Network Design Problem} \label{subsec:prob_formulation}

The goal of the maximum accessibility network design problem (MaxANDP) is to design a cycling network subject to a fixed budget such that the total accessibility of a given set of OD pairs, denoted by $\mS$, is maximized. Such a set may be defined based on geographical units \citep{imani2019cycle} or ridership data \citep{liu2021urban}. Various metrics have been proposed to measure accessibility, mostly focusing on first finding one or more routes between each OD pair using the designed network and then calculating the accessibility based on the selected routes. 

Let $\mG = (\mN, \mE)$ be a directed graph where $\mE$ is the set of edges and $\mN$ is the set of nodes, corresponding to road segments and intersections, respectively. Each edge $(i, j)\in \mE$ is assigned a travel time $t_{ij}$. We denote by $\mE^+(i)$ and $\mE^-(i)$ the sets of incoming and outgoing edges of node $i$, respectively. Edges and nodes are partitioned into high-stress and low-stress sets according to a cycling stress assessment based on road geometry, existing infrastructure, and vehicle traffic conditions \citep{furth2016network}. We assume that cyclists prefer cycling on low-stress roads over high-stress roads. Sets with subscripts $h$ and $l$ indicate the high-stress and low-stress subsets of the original set, respectively. High-stress edges $(i, j)\in \mE_{h}$ and nodes $i\in \mN_{h}$ are assigned costs $c_{ij}$ and $b_{i}$, respectively, representing the costs of turning them into low-stress through building new infrastructure, e.g., cycle tracks or traffic lights.

Let $\bfx \in \{0, 1\}^{|\mE_h|}$ and $\bfz\in \{0, 1\}^{|\mN_h|}$, respectively, denote the \textit{edge selection} and \textit{node selection} variables (referred to as \textit{network design} decisions), whose components are 1 if that edge or node is chosen for the installation of infrastructure that makes it low stress. Edge and node selections are subject to budgets $B_\text{edge}$ and $B_\text{node}$, respectively. Let $\bfy^{od}\in \{0, 1\}^{|\mE|}$ denote the \textit{routing decision} associated with OD pair $(o, d) \in \mS$. The routing problem on a network specified by $\bfx$ and $\bfz$ is characterized by an objective function $h^{od}(\bfx, \bfz, \cdot): \{0, 1\}^{|\mE|} \rightarrow \bbR$ and a feasible set $\mY^{od}(\bfx, \bfz) \subseteq \{0, 1\}^{|\mE|}$. A function $g(\bfx, \bfz, \cdot): \{0, 1\}^{|\mE|} \rightarrow \bbR_+$ is used to calculate the accessibility of each OD pair based on the selected route(s). Each OD pair is weighted by a constant $q^{od}\in \bbR_+$ (e.g., population). The MaxANDP is formulated as
\begin{subequations}
\label{prob:general_cycling}
\begin{align}
    \underset{\bfx, \bfz, \bfy^{od}}{\textrm{maximize}} \quad
        & \sum_{(o, d) \in \mS} q^{od} g(\bfx, \bfz, \bfy^{od}) 
        \label{cyc_general:obj} \\
    \textrm{subject to} \quad
        & \bfy^{od} \in 
            \argmin_{\bfy \in \mY^{od}(\bfx, \bfz)} h^{od}(\bfx, \bfz, \bfy),
            \quad (o, d) \in \mS
        \label{cyc_general:opt_route} \\
        \label{cyc_general:edge_budget}
        & \bfc^\intercal \bfx \leq B_\text{edge} \\
        \label{cyc_general:node_budget}
        & \bfb^\intercal \bfz \leq B_\text{node} \\
        \label{cyc_general:domain}
        & \bfx\in \{0, 1\}^{|\mE_h|}, 
        \bfz\in \{0, 1\}^{|\mN_h|},
\end{align}
\end{subequations}
where $\bfc$ and $\bfb$ indicate cost vectors for high-stress edges and nodes, respectively. The objective function~\eqref{cyc_general:obj} maximizes total cycling accessibility. Constraints~\eqref{cyc_general:opt_route} ensure that the selected routes are optimal for the OD pairs' objective functions. Constraints~\eqref{cyc_general:edge_budget} and \eqref{cyc_general:node_budget} enforce budgets on the network design. We emphasize that followers in Problem~\eqref{prob:general_cycling} are independent of each other, in contrast to in other transportation network design problems where the followers typically interact with each other via traffic equilibrium conditions. This is a common practice in cycling network design because i) unlike motor traffic, cycling traffic congestion is less common in practice, and ii) the goal of Problem~\eqref{prob:general_cycling} is not to model travel behavior, but to quantify network accessibility, i.e., how far a cyclist can reach in a specified time. The objective of Problem~\eqref{prob:general_cycling} has been shown to be highly correlated with the cycling mode choice in the City of Toronto \citep{imani2019cycle}, and is being used as an evaluation metric in support of cycling infrastructure project prioritization \citep{city2021active}.

Many variants of Problem~\eqref{prob:general_cycling} have been proposed in the literature, corresponding to different combinations of accessibility measure (specified by $g$) and routing problem (specified by $h^{od}$ and $\mY^{od}$). To illustrate our method, we consider two problems: i) one that uses location-based accessibility measures and shortest-path routing problems, and ii) one proposed by \citet{liu2021urban} that employs a utility-based accessibility measure and discrete route choice models. These two problems allow us to illustrate the generality of our approach in terms of the follower cost function: the former problem has follower cost $g(\bfy)$ and the latter uses $g(\bfx, \bfz, \bfy)$, that is, the leader decision is also included in the cost. The former problem is locally more relevant in Toronto and we briefly describe it next. We refer readers to \citet{liu2021urban} for more details on the latter problem.

Location-based accessibility measures use a decreasing function of the travel time from origin to destination, namely an impedance function, to model the dampening effect of separation. We consider a piecewise linear impedance function
\begin{equation}
\label{eq:loc_acc}
    g(\bfy^{od}) = 
    \begin{cases}
       1 - \beta_1 \bft^\intercal \bfy^{od}, 
            & \textrm{ if $\bft^\intercal \bfy^{od} \in [0, T_1)$}\\
        1 - \beta_1 T_1 - \beta_2 \left(\bft^\intercal \bfy^{od} - T_1\right), 
           &  \textrm{ if $\bft^\intercal \bfy^{od} \in [T_1, T_2)$}\\
        0, 
            & \textrm{ if $\bft^\intercal \bfy^{od} \geq T_2$},
    \end{cases}
\end{equation}
where $\bft$ indicates a vector of edge travel times, $T_1, T_2\in \bbR_+$ are breakpoints, and $\beta_1, \beta_2 \in \bbR_+$ are penalty factors for intervals $[0, T_1)$ and $[T_1, T_2)$, respectively. 
This function can be used to approximate commonly used impedance functions, including negative exponential, rectangular, and linear functions (visualized in \ref{appsub:acc_measure}). 
While we consider two breakpoints for simplicity, the formulation can be easily generalized to account for more. 


We use the level of traffic stress (LTS) metric \citep{furth2016network} to formulate the routing problems.
Let $\mathbf{A}$ be the node-edge matrix describing the flow-balance constraints on $\mG$, and $\bfe^{od}$ be a vector whose $o^\text{th}$ and $d^\text{th}$ entries are 1 and $-1$, respectively, with all other entries 0. Given network design $(\bfx, \bfz)$, the routing problem for $(o, d) \in \mS$ is formulated as
\begin{subequations}
\label{prob:od_routing}
\begin{align}
\underset{\bfy^{od} \in \{0, 1\}^{|\mE|}}{\textrm{minimize}} \quad
    \label{od_routing:obj}
    & \bft^\intercal \bfy^{od} \\
\textrm{subject to} \quad 
    \label{od_routing:balance}
    & \mathbf{A}\bfy^{od} = \bfe^{od} \\
    \label{od_routing:edge_design}
    & y^{od}_{ij} \leq x_{ij},
        \quad \forall (i,j)\in \mE_{h}\\
    \label{od_routing:node_design}
    & y^{od}_{ij} \leq x_{wl} + z_i, 
        \quad \forall i \in \mN_h, 
        \ (i,j)\in \mE^-_{h}(i), 
        \ (w,l)\in \mE^-_{h}(i) \cup \mE^+_{h}(i).
\end{align}
\end{subequations}

Objective function~\eqref{od_routing:obj} minimizes the travel time. Constraints~\eqref{od_routing:balance} direct one unit of flow from $o$ to $d$. Constraints~\eqref{od_routing:edge_design} ensure that a currently high-stress edge can be used only if it is selected. Constraints~\eqref{od_routing:node_design} guarantee that a currently high-stress node can be crossed only if either the node is selected or all high-stress edges that are connected to this node are selected. This is an exact representation of the intersection LTS calculation scheme that assigns the low-stress label to a node if traffic signals are installed or all incident roads are low-stress \citep{imani2019cycle}. To ensure Problem~\eqref{prob:od_routing} is feasible, we add a virtual low-stress link from $o$ to $d$ and set its travel time to $T_2$. In doing so, the travel time is $T_2$ when the destination is unreachable on the low-stress network, corresponding to zero accessibility, as defined in equation~\eqref{eq:loc_acc}. The full formulation is in \ref{app:maxandp_formulation}.

Solving Problem~\eqref{prob:general_cycling} is challenging due to the large number of discrete decisions in the leader's problem, which is further compounded by a large $\mS$. When using the location-based accessibility measures defined in \eqref{eq:loc_acc}, we can adapt the the Benders decomposition approach from \citet{magnanti1986tailoring} to solve it (detailed in \ref{app:benders}). When using the utility based measure, we can apply the algorithm from \citet{liu2021urban}. Although these algorithms allow us to solve the original bilevel model to optimality for synthetic instances (Section \ref{subsec:exp_opt}), they are insufficient for real instances with over one million followers (Section \ref{sec:case_study}), as the problem size grows quickly with $|\mS|$, necessitating the use of sampling to reduce the problem size. \revision{Since Problem~\eqref{prob:pwo_ext} shares the same structure as the general bilevel model~\eqref{prob:general}, we can use the reduced model~\eqref{prob:reduced} and our ML-augmented model~\eqref{prob:pwo_ext} to generate leader solutions for it.}

\subsection{Learning Follower Representations} \label{sec:repre_learning}

\revision{Applying our ML-augmented model requires using follower features as inputs to the ML model. 
In this section, we introduce a learning framework that maps the set of followers $\mS$ to a $\xi$-dimensional feature space. We present our approach in Section \ref{subsubsec:our_fea} and discuss its alternatives in Section \ref{subsubsec:alt_fea}.}

\subsubsection{An NLP-inspired approach.} \label{subsubsec:our_fea}

\revision{The design of our follower embedding method is informed by our theoretical analysis. Recall that the bounds in Theorems \ref{thm:non-parametric_bound} and \ref{thm:reg_bound} depend on the Lipschitz constant $\mu$, which bounds the rate at which the follower cost changes as a function of changes in the feature space. To reduce $\mu$ and thus tighten these bounds, it is essential to construct a feature space where ``similar followers'' are mapped close to each other. This insight draws an analogy to the Natural Language Processing (NLP) literature, where the problem of ``pulling similar words together'' has been extensively studied \citep{mikolov2013distributed, radford2018improving}. Next, we introduce a two-step framework that transforms the follower representation problem into a word embedding problem, allowing us to leverage well-established techniques in the NLP literature.}

\begin{itemize}
    \item[\textbf{Step I:}] \textbf{Relationship graph construction}. We begin by constructing a relationship graph $\mR$, where each node represents a follower, and each edge is weighted to reflect the similarity between followers. Many metrics have been proposed to quantify such similarities, with most focusing on the ``opportunity cost'' of switching between followers. For example, \citet{keutchayan2021problem} define the opportunity cost of applying the leader decision that is optimal to scenario $t$ in scenario $s$ as $d(s, t) = G^s(\bfx^{t*}) - G^s(\bfx^{s*})$ where $\bfx^{s*}$ and $\bfx^{t*}$ denote the optimal leader solutions obtained by solving the single-follower version of Problem~\eqref{prob:general_cycling} with followers $s$ and $t$, respectively. Building on a similar idea, \citet{bertsimas2022optimization} define a symmetric metric for followers $s$ and $t$ as $\left(d(s, t) + d(t, s)\right)/2$. While these metrics have been shown to be effective and are compatible with our framework, they require solving $|\mS|$ single-scenario versions of Problem~\eqref{prob:general}, which is computationally expensive when $\mS$ is huge and when the leader's problem is nonconvex.

Motivated by the fact that evaluating follower costs given a leader's solution is computationally cheaper, we propose a new approach that quantifies follower similarity based on their costs under some sampled leader solutions. Specifically, we randomly sample $n_\text{sim}$ leader decisions $\{\bfx_i \in \mX \}_{i=1}^{n_\text{sim}}$. For each sampled $\bfx_i$, we calculate the costs $G^s_i := G^s(\bfx_i)$ for all followers $s\in \mS$ by solving their follower problems. We then define the weight of the edge between followers $s, t\in \mS$ as
\begin{equation} \label{eq:weight}
    \pi^{st} := 
    \exp \left[ 
        \frac{1}{n_\textrm{sim}}
        \sum_{i = 1}^{n_\textrm{sim}} 
        - |G^s_i - G^t_i|
    \right].
\end{equation}

    \item[\textbf{Step II:}] \textbf{Follower embedding}. Once $\mR$ is constructed, we adapt the DeepWalk algorithm proposed by \citet{perozzi2014deepwalk} to learn node representations. We first generate a set of random walks in $\mR$, and then apply the SkipGram algorithm \citep{mikolov2013distributed}, which was designed for learning word embeddings, to learn node features treating each node and each random walk as a word and a sentence, respectively. Unlike \citet{perozzi2014deepwalk}, who generate random walks by uniformly sampling nodes connected to the current node, we generate random walks according to the weights assigned to edges incident to the current node. So, followers that yield similar results under the sampled leader decisions are likely to appear in a same walk, and thus will be close to each other in the feature space.

\end{itemize}

\subsubsection{Alternative approaches.} \label{subsubsec:alt_fea} 
\revision{To the best of our knowledge, the problem of learning follower features has not been studied in the literature. However, several ideas from two-stage stochastic programming, continuous approximation, and graph theory may be adapted to our context. Below, we discuss these ideas, two of which serve as baselines in Section \ref{subsec:exp_pred}.}

\begin{enumerate}
    \item \revision{\textbf{Scenario reduction.} A common approach in scenario reduction is to represent each scenario using the realized second-stage parameters. However, in Problem~\eqref{prob:general_cycling}, this representation corresponds to the right-hand side of the flow-balance constraints~\eqref{od_routing:balance}, whose dimensionality equals the number of nodes in the road network and whose entries are mostly zeros. Such a feature representation is high-dimensional and uninformative, making ML-augmented models intractable and rendering the sampling method ineffective (recall discussions in Sections \ref{subsubsec:ml_para_reg} and \ref{subsubsec:prac_non-parametric}).}

    \item \revision{\textbf{TSP features.} In Problem~\eqref{prob:general_cycling}, our ML prediction target is a function of the optimal value of a shortest path Problem~\eqref{prob:od_routing}. This problem can be viewed as a special case of the traveling salesman problem (TSP), where the number of stops is set to two. We thus can adapt the travel time predictors proposed by \citet{liu2021time}, which are well grounded in the continuous approximation literature for predicting TSP objective values. These predictors include the locations of the origin and destination, as well as geometric features describing their relative positions (e.g., the distance between them), as detailed in \ref{ec:tsp_fea}. The TSP features are included as a baseline.}

    \item \revision{\textbf{Graph-theoretical (GT) features.} An alternative approach to leveraging the graph structure established in Step I of our method is to represent each follower using node importance metrics from graph theory. To demonstrate the effectiveness of our DeepWalk algorithm, we include a set of these metrics as baseline features in our computational studies. These features, including node centrality, betweenness, and their variants, are detailed in \ref{ec:graph_fea}.}
\end{enumerate}

\subsection{Experiment 1: Predicting OD-Pair Accessibility Using ML Models}
\label{subsec:exp_pred}
In this section, we conduct experiments to validate (i) the effectiveness of our follower sampling method in improving prediction accuracy and (ii) the predictive power of our learned features.

\textbf{Experiment setup.} We consider four accessibility measures as prediction targets: three location-based measures using exponential (EXP), linear (LIN), and rectangular (REC) impedance functions, and one utility-based (UT) measure from \citet{liu2021urban}. For each measure, we randomly generate 3,000 network designs and calculate the accessibility of each follower under every design. The accessibility associated with each network design constitutes a dataset, which we split into training and testing sets to train ML models and evaluate their prediction performance. The evaluation includes three ML models compatible with our ML-augmented framework: $k$NN, lasso regression, and ridge regression. We use the mean absolute error (MAE), normalized by the average total accessibility across the 3,000 network designs, as the evaluation metric. \revision{We vary the training sample size between 1\%--5\% of all OD pairs because implementation of both the reduced and ML-augmented models on large real-world case studies (Section \ref{sec:case_study}) are only possible when the sample size is very small ($<0.2\%$ in our case study)}.

\textbf{Baselines.} For sampling, we evaluate our balanced $p$-median sampling method (BMED) and compare it against four commonly used baselines from the scenario reduction literature: uniform sampling (UNI), $p$-center sampling (CEN), \revision{Euclidean Wasserstein scenario reduction (EW), and the method proposed by \citet{dupavcova2003scenario} (DPCV)}. 
Since the BMED, CEN, and EW problems are $\mathcal{NP}$-hard, we adapt heuristics from \citet{boutilier2020ambulance}, \citet{gonzalez1985clustering}, and \citet{bertsimas2022optimization} to solve them, respectively. 
These methods involve randomness, so we repeat each approach 10 times with different random seeds and report the mean and 95\% confidence intervals of the normalized MAE. 
\revision{For features, we benchmark our representation learning-based features (REP) against the TSP and GT features introduced in Section \ref{subsubsec:alt_fea}.}



\begin{figure}[!ht]
    \centering
    \includegraphics[width=0.75 \textwidth]{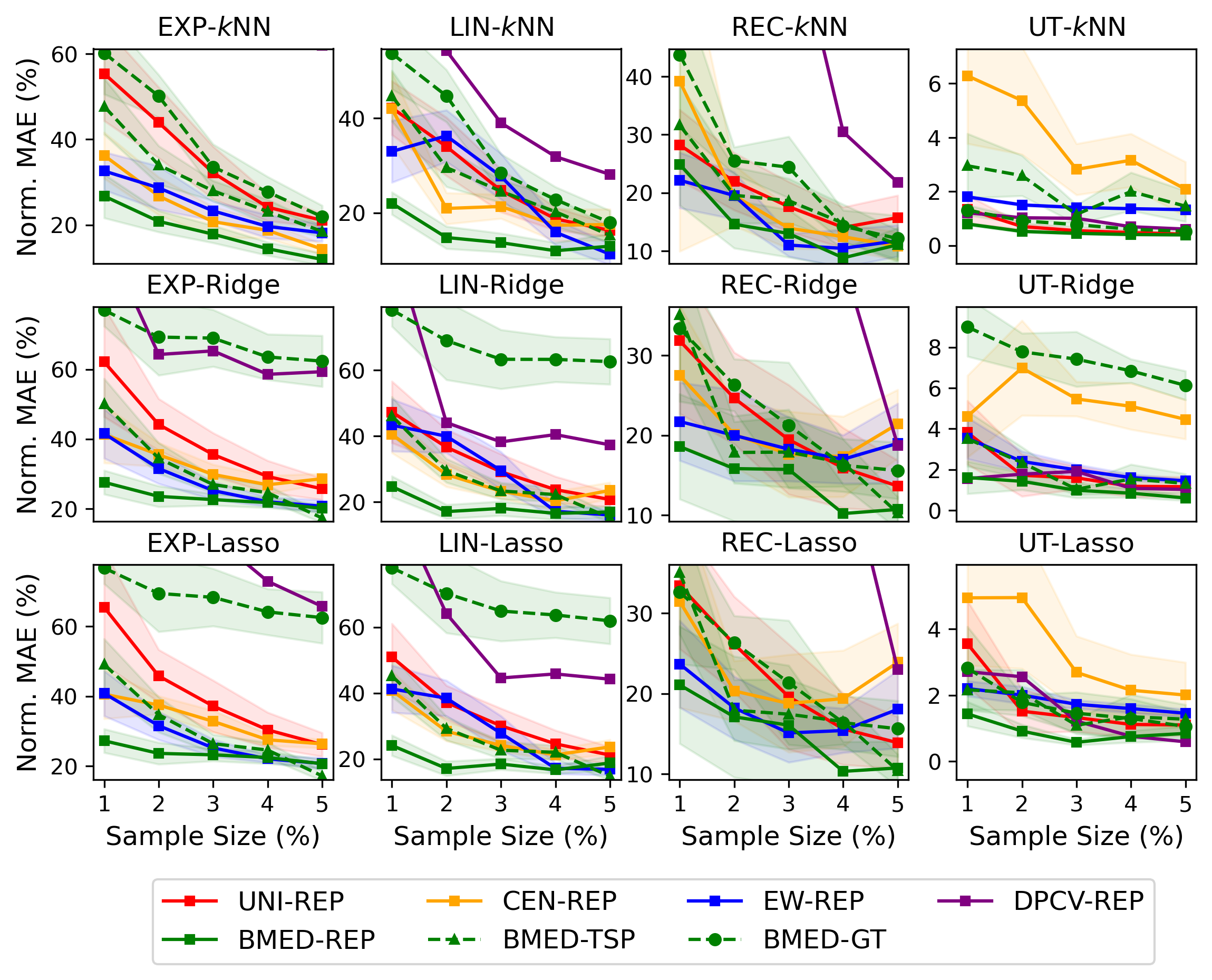}
    \caption{Normalized MAE ($\pm$ 95\% confidence interval) over 3,000 network designs. Each panel corresponds to an ``accessibility measure'' - ``ML model'' pair. Each line indicates a ``sampling method'' - ``feature'' combination.
    For readability, we only present the best-performing sampling method for TSP and GT features. Some lines are partially visible as they are significantly inferior to others.}
    \label{fig:online_pred}
\end{figure}

\textbf{Effectiveness of the follower selection algorithm.} As illustrated in Figure \ref{fig:online_pred}, our sampling method, BMED, consistently achieves the lowest prediction error when using REP features, regardless of the accessibility measures, ML models, or sample sizes. The performance advantage is particularly pronounced with extremely small sample sizes (e.g., 1\%), where BMED outperforms baseline strategies by over 20\%, highlighting the effectiveness of our bounds in guiding sample selection. Furthermore, BMED generally exhibits lower variation in prediction error compared to the baselines. Similar trends are observed for TSP and GT features.

\textbf{Predictive power of the learned features.} ML models generally perform better with our REP features compared to baselines. As shown in Figure \ref{fig:online_pred}, the REP features significantly outperform the best baseline (TSP) when using $k$NN, with improvements of 44.0\%, 50.7\%, 21.1\%, and 73.0\% for the four accessibility measures, respectively, at a 1\% sample size. 
With lasso and ridge regression, REP features also outperform TSP features, underscoring the robustness of our representation learning approach. 
\revision{We note that the performance gap between BMED-REP and BMED-TSP decreases as the sample size increases. This observation aligns with empirical evidence from the representation learning literature, where the benefit of ``high-quality'' features is usually more pronounced when the sample size is small \citep{perozzi2014deepwalk, cole2022does}.}
\revision{Finally, the prediction error achieved by the REP features is typically 75\% lower than that of GT features, demonstrating the effectiveness of our DeepWalk algorithm in leveraging the follower relationship graph.}


\subsection{Experiment 2: Generating Leader Decisions using ML-augmented Models} \label{subsec:exp_opt}

Next, we evaluate the ability of our ML-augmented model, enhanced by our follower sampling and representation learning algorithms, to generate high-quality leader solutions.

\textbf{Experiment setup.} We create 12 problem instances on a synthetic network, corresponding to combinations of three design budgets and four accessibility measures. To compute the optimality gap of the derived leader solutions, we adapt the Benders approach from \citet{magnanti1986tailoring} to solve the synthetic instances to optimality (details in \ref{app:benders}). Each model is applied 10 times with samples generated using different random seeds, and we report the average optimality gap. 

\textbf{Baselines.} We implement our ML-augmented model with $k$NN and linear regression and benchmark them against the reduced model~\eqref{prob:reduced} using UNI, BMED, CEN, EW, and DPCV samples.


\begin{figure}[!ht]
    \centering
    \includegraphics[width=.8\textwidth]{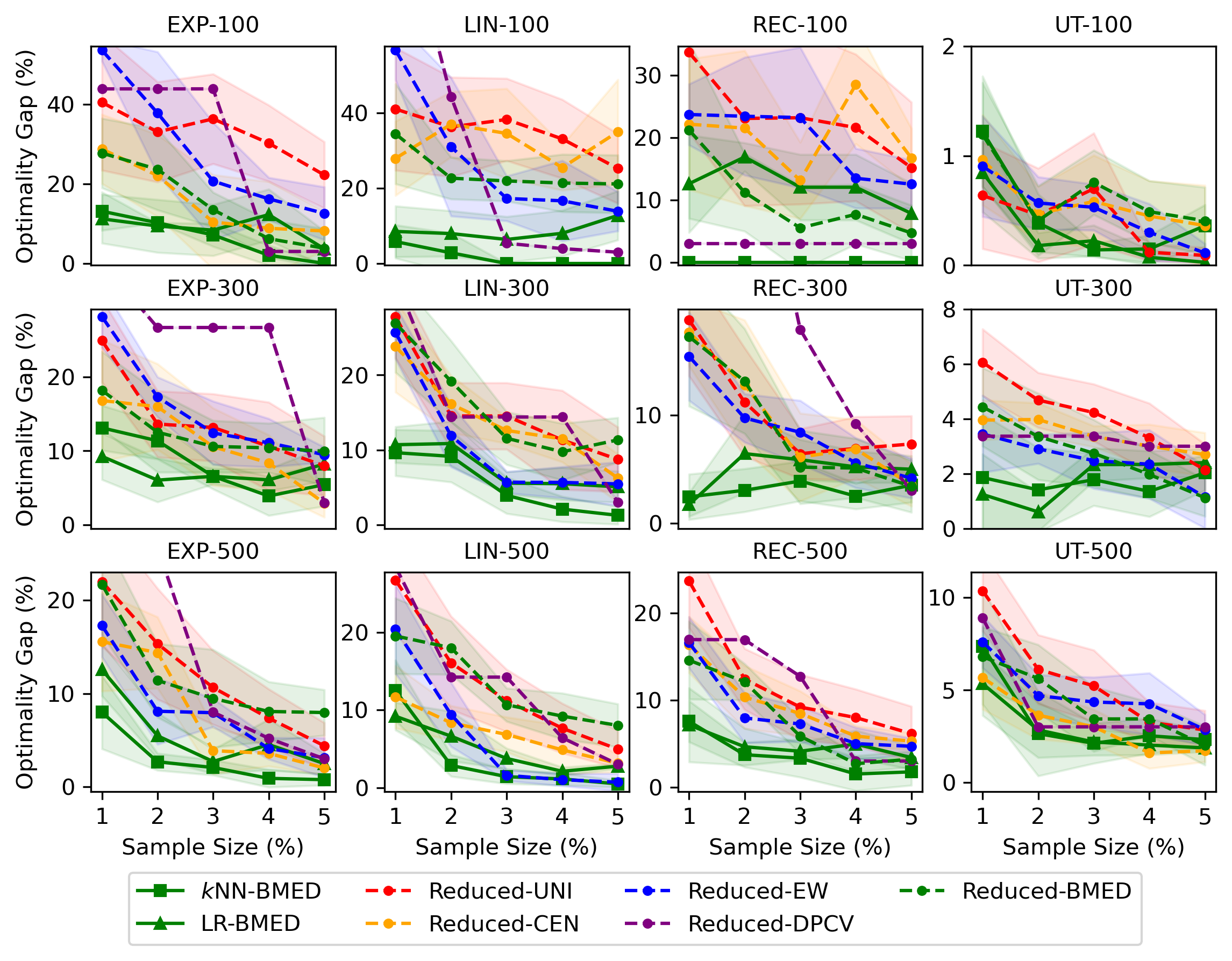}
    \caption{Mean optimality gap ($\pm$ 95\% confidence interval) of leader solutions. Problem instances are named as ``accessibility measure''-``budget'' and the solution methods are named as ``model''-``sampling method''. Optimality gap $= |F(\bfx^*) - F(\bfx')| / F(\bfx^*)$ where $\bfx^*$ and $\bfx'$ are leader decisions from Benders decomposition and a sampling-based method, respectively. Some lines are partially visible as they are significantly inferior to other others.}
    \label{fig:OptGap}
\end{figure}

\textbf{Effectiveness of the ML-augmented models.} As shown in Figure \ref{fig:OptGap}, ML-augmented models ($k$NN-BMED and REG-BMED) generally outperform the reduced models by a large margin, especially when the sample size is extremely small (e.g., 1\%). Additionally, the confidence intervals of the ML-augmented models are consistently narrower than those of the reduced models, indicating higher stability. The ML component effectively captures the impact of leader decisions on unsampled followers, resulting in solutions of higher quality and robustness.



\begin{figure}
    \centering
    \includegraphics[width=.8\textwidth]{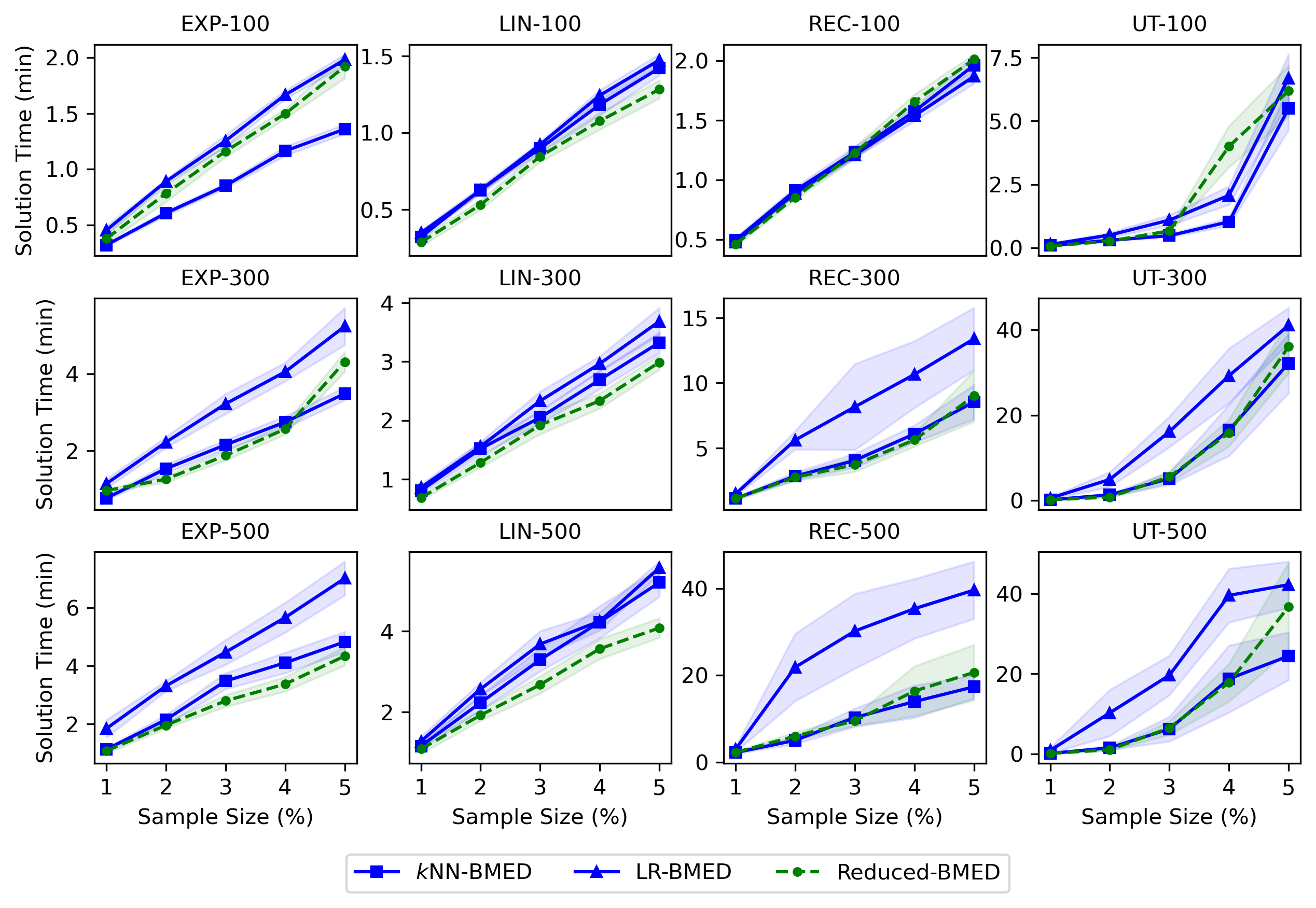}
    \caption{Mean solution time ($\pm$ 95\% confidence interval). Problem instances are named as ``accessibility measure''-``budget'' and the solution methods are named as ``model''-``sampling method''.}
    \label{fig:sol_time}
\end{figure}

\textbf{Efficiency of the ML-augmented models.} Figure \ref{fig:sol_time} presents the solution time of the three models with BMED samples. In general, the solution time of all models increases as the sample size increases. The $k$NN-augmented model and the reduced model require similar solution time as the former is a re-weighted version of the latter and does not have any additional decision variables or constraints. The linear regression-augmented model requires longer solution time as it has more decision variables. Compared to applying Benders decomposition to the original model which generally takes over 10 hours for each instance, the ML-augmented models generate leader decisions of similar quality in 0.5--5\% of the solution time, highlighting the efficiency of our method.

\section{Case Study: Cycling Infrastructure Planning in the City of Toronto} \label{sec:case_study}
In this section, we present a case study applying our method to Toronto, Canada. 
We started a collaboration with the City's Transportation Services Team in September 2020, focusing on developing quantitative tools to support cycling infrastructure planning in Toronto. As an evaluation metric, low-stress cycling accessibility has been used by the City of Toronto to support project prioritization \citep{city2021plan}. Our optimization model, which maximizes low-stress cycling accessibility, was used to support Toronto's 2025--2030 cycling infrastructure planning \citep{CityofToronto2024b}.
We introduce Toronto's cycling network in Section \ref{subsec:trt_cyc_network} and use our methodology to examine actual and future decisions regarding network expansion in Section \ref{subsec:exp_trt}. 

\subsection{Cycling Network in Toronto} \label{subsec:trt_cyc_network}
\begin{figure}[!ht]
    \centering
    \includegraphics[width=.7\textwidth]{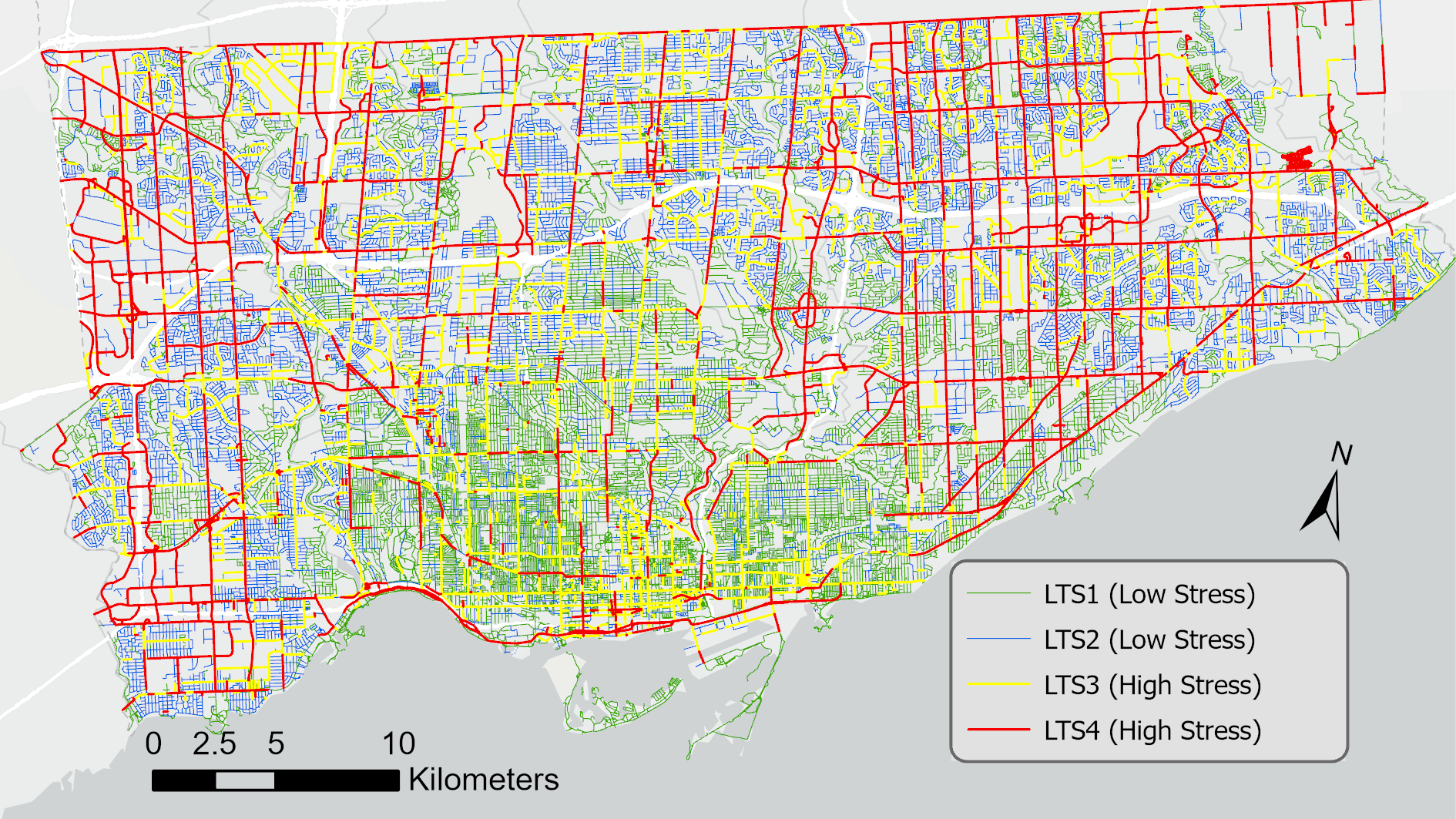}
    \caption{Level of traffic stress of Toronto's road network (July 2021).}
    \label{fig:lts_network}
\end{figure}

We construct Toronto's cycling network based on the centerline network retrieved from the Toronto Open Data Portal \citep{trt2020open}. We pre-process the network by removing roads where cycling is legally prohibited, deleting redundant nodes and edges, and grouping arterial roads into candidate cycling infrastructure projects (detailed in \ref{appsub:preprocessing}). The final cycling network has 10,448 nodes, 35,686 edges, and 1,296 candidate projects totaling 1,913 km. We use the methods summarized in \citet{lin2021impact} to calculate the LTS of each link in the cycling network. LTS1 and LTS2 links are classified as low-stress, while LTS3 and LTS4 links are high-stress since LTS2 corresponds to the cycling stress tolerance for the majority of the adult population  \citep{furth2016network}. 
Although most local roads are low-stress, high-stress arterials create many disconnected low-stress ``islands'', limiting low-stress cycling accessibility in many parts of Toronto (see Figure \ref{fig:lts_network}). 

We use the following procedure to calculate the low-stress cycling accessibility of Toronto, which serves as an evaluation metric of Toronto's cycling network and the objective of our cycling network design Problem~\eqref{prob:general_cycling}. The city is divided into 3,702 geographical units called dissemination areas (DAs). We define each DA centroid as an origin with every other DA centroid that is reachable within 30 minutes on the overall network being a potential destination, totaling 1,154,663 OD pairs ($\mS$). These OD pairs are weighted by the job counts at the destination ($q^{od}$), retrieved from the 2016 Canadian census \citep{StatsCan2016}. We use a rectangular impedance function with a cut-off time of 30 minutes ($g$). We assume a constant cycling speed of 15 km/h for travel time calculation. The resulting accessibility measure can be interpreted as the total number of jobs that one can access within 30 minutes via low-stress routes in the City of Toronto. This metric has been shown to be highly correlated with cycling mode choice in Toronto \citep{imani2019cycle}.

\subsection{Expanding Toronto's Cycling Network} \label{subsec:exp_trt}
\begin{figure}[!ht]
    \centering
    \includegraphics[width=.7\textwidth]{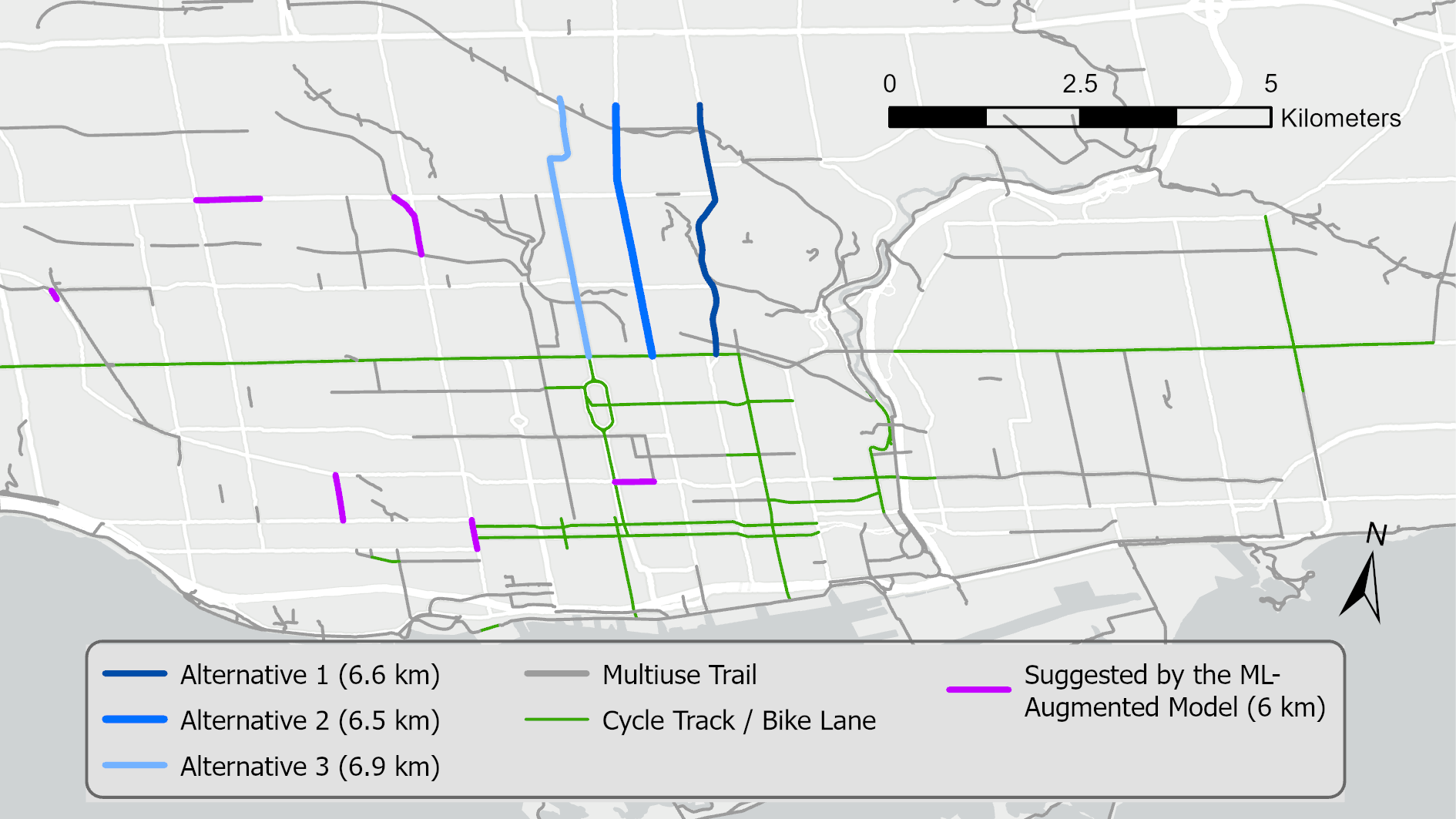}
    \caption{Project alternatives and the existing cycling infrastructure in the City of Toronto (January 2021).}
    \label{fig:trt_alternatives}
\end{figure}

As a part of our collaboration, in January 2021 we were asked to evaluate the accessibility impact of three project alternatives for building bike lanes (see Figure \ref{fig:trt_alternatives}) to meet the direction of Toronto’s City Council, intended to provide a cycling connection between midtown and the downtown core \citep{city2021active}. These projects were proposed in 2019 but their evaluation and implementation were accelerated because of increased cycling demand during COVID. We determined that alternative 2 had the largest accessibility impact. It was  
ultimately implemented due to its accessibility impact and other performance indicators \citep{city2021active}.

This decision-making process exemplifies the current practice of cycling infrastructure planning in Toronto: i) manually compile a list of candidate projects, ii) rank the candidate projects based on certain metrics, and iii) design project delivery plans. From a computational perspective, steps i) and ii) serve as a heuristic for solving MaxANDP. This heuristic approach was necessary for several reasons, including political buy-in for the three alternatives, and the computational intractability of solving MaxANDP at the city scale. In fact, Benders decomposition, which was used to solve the synthetic instances in Section \ref{sec:comp}, cannot find a feasible solution to these instances before running out of memory. Now, we can use our ML-augmented model to search for project combinations at the city scale without pre-specifying project candidates.

To this end, we first apply the ML-augmented model with a budget of 6 km (similar to alternative 2). The optimal projects (see Figure \ref{fig:trt_alternatives}) improve Toronto's total low-stress cycling accessibility by 9.46\% over alternative 2. Instead of constructing only one corridor as in alternative 2, the ML-augmented model selects six disconnected road segments. Some of them serve as connections between existing cycling infrastructure, others bridge currently disconnected low-stress sub-networks consisting of low-stress local roads. We also compare our approach against i) three reduced models and ii) a greedy heuristic that iteratively selects the candidate project that leads to the maximum increase in total accessibility until the budget is depleted. As presented in Table \ref{tab:real_alg_comp}, the greedy heuristic, which is commonly adopted in practice and in the existing literature, closely matches the performance of the human-proposed solution \revision{(i.e., alternative 2 identified through the process described above)}. With similar computational times, the three reduced models are all inferior to our model, with the best reduced model being on par with the human performance and others lagging over 20\% behind. Interestingly, the greedy heuristic performs quite well against the reduced model. We believe this highlights the difficulty of achieving strong performance with a small sample in a purely sampling based model.

\begin{table}[!ht]
\centering
\small
\caption{Increases in the average low-stress cycling accessibility over 3,702 DAs in Toronto due to 6 km of new cycling infrastructure selected by different approaches. We implement each sampling-based method with five random samples and report the best result across the five samples.} \label{tab:real_alg_comp}
\begin{tabular}{@{}lrr@{}}
\toprule
\multicolumn{1}{l}{Method (Sample)} & \multicolumn{1}{r}{Accessibility Increase} &  \multicolumn{1}{r}{\% change relative to human} \\ \midrule
Human                      &  6,902 &    $+\ 0.00\%$ \\
Greedy                     &  7,012 &    $+\ 1.59\%$      \\
Reduced (UNI)              &  4,965 &    $-\ 28.06\%$ \\
Reduced (PCEN)             &  5,730 &    $-\ 20.45\%$ \\
Reduced (BMED)             &  7,118 &    $+\ 3.13\%$  \\ 
\textbf{ML-augmented (BMED)} &  \textbf{7,555} & $\mathbf{+9.46\%}$\\ \bottomrule
\end{tabular}
\end{table}

Next, to demonstrate the potential impact of our method in Toronto, we compare our model versus the greedy heuristic, which mimics the current cycling infrastructure planning practice. We increase the road design budget from 10 to 100 km in increments of 10 km. The 100 km budget aligns with Toronto's cycling network expansion plan for 2022--2024 \citep{city2021plan}. The greedy heuristic took over 3 days to expand the network by 100 km as each iteration involves solving millions of shortest path problems. Our approach took around 4 hours to find a leader decision using a sample of 2,000 OD pairs (1.7\% of all OD pairs). Given this speedup, we can solve our model multiple times with different samples and report the best solution as measured by the total accessibility of all OD pairs. The computational setups of the greedy heuristic and our approach are detailed in \ref{appsub:trt_comp}. 

\begin{figure}[!h]
    \centering
    \includegraphics[width=2in]{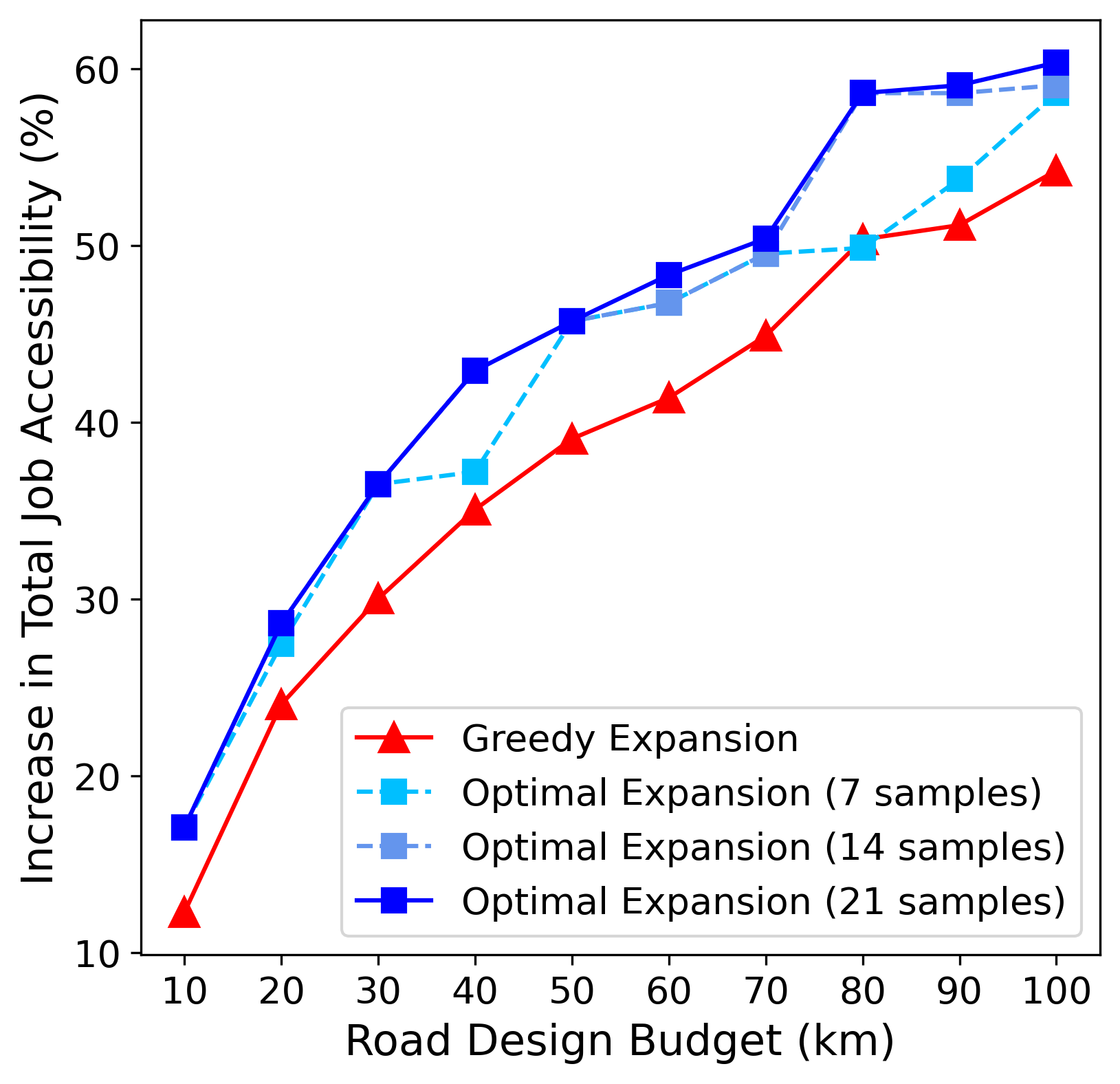}
    \caption{The performance profiles of the greedy and optimal expansions. Note that using 21 different sets of OD pairs samples results in the same solution time as the greedy expansion. Hence, 7 and 14 samples correspond to $1/3$ and $2/3$ of the solution time of greedy.}
    \label{fig:opt_vs_greedy}
\end{figure}

As shown in Figure \ref{fig:opt_vs_greedy}, when holding both methods to the same computational time (meaning that we solve our ML-augmented model with 21 different sets of OD pair samples and taking the best solution), our approach increases accessibility by 19.2\% on average across different budgets. For example, with a budget of 70 km, we can improve the total accessibility by a similar margin as achieved by the greedy heuristic using a 100 km budget, corresponding to a savings of 18 million Canadian dollars estimated based on the City's proposed budget \citep{city2021plan}. If instead we used the full 100 km budget, we would achieve 11.3\% greater accessibility. Finally, we note that solution quality was similar between 14 and 21 samples, meaning that with we can achieve the above gains while simultaneously reducing solution time by approximately 33\%.

\subsection{Managerial Insights}
\revision{We conclude the case study by providing managerial insights derived by examining the network design solutions generated by our approach and the current practice (visualized in \ref{appsub:greedy_opt_comp}), and through conversations with practitioners. }

\revision{First, both methods prioritize cycling infrastructure projects in the downtown area, where a well-connected low-stress cycling network already exists (Figure \ref{fig:lts_network}) and where job opportunities are densely distributed. These projects connect many DAs to the existing low-stress network, granting them access to a large number of job opportunities nearby. This empirical finding underscores the importance of forming a well-connected cycling network and aligns with the City of Toronto's emphasis on enhancing cycling network connectivity \citep{city2021active, CityofToronto2024}.}

\revision{Second, the improvements achieved by our approach mainly come from identifying projects that have little accessibility impact when constructed alone, yet exhibit significant impact when combined. These improvements often involve eliminating high-stress barriers between low-stress cycling ``islands'', which requires coordinated efforts across multiple infrastructure projects. This observation highlights the potential for improving the City's infrastructure investment outcomes through city-scale joint planning. While such coordination was previously impossible due to the computational challenges discussed earlier, our machine learning-augmented optimization method provides a scalable approach to identifying such solutions going forward.}

Finally, although maximizing low-stress cycling accessibility is not the only goal of cycling network design, we believe our approach can be useful in at least three contexts: i) in the long term, our model can be used to generate a base plan that can later be fine-tuned by transportation planners; ii) in the near term, our approach can efficiently search for project combinations from a large pool that would be very difficult to analyze manually; iii) Given a fixed budget, our model provides a strong benchmark against which to validate the goodness of human-proposed solutions. Due to its strong performance, our optimization model was used to inform Toronto's 2025--2030 cycling infrastructure plan \citep{CityofToronto2024, CityofToronto2024b}.

\section{Conclusion}
In this paper, we present a novel ML-based approach to solving bilevel (stochastic) programs with a large number of independent followers (scenarios). We build on two existing strategies---sampling and approximation---to tackle the computational challenges imposed by a large follower set. The model considers a sampled subset of followers while integrating an ML model to estimate the impact of leader decisions on unsampled followers. Unlike existing approaches for integrating optimization and ML models, we embed the ML model training into the optimization model, which allows us to employ general follower features that may not be compactly represented by leader decisions. Under certain assumptions, the generated leader decisions enjoy solution quality guarantees as measured by the original objective function considering the full follower set. We also introduce practical strategies, including follower sampling algorithms and a representation learning framework, to enhance the model performance. Using both synthetic and real-world instances of a cycling network design problem, we demonstrate the strong computational performance of our approach in generating high-quality leader decisions. The performance gap between our approach and baseline approaches are particularly large when the sample size is small. 


\ACKNOWLEDGMENT{The authors are grateful to Sheng Liu, Merve Bodur, Elias Khalil, Rafid Mahmood, and Erick Delage for helpful comments and discussions. This research is supported by funding from the City of Toronto and NSERC Alliance Grant 561212-20. Resources used in preparing this research were provided, in part, by the Province of Ontario, the Government of Canada through CIFAR, and companies sponsoring the Vector Institute.}


\bibliographystyle{informs2014} 
\bibliography{paper.bbl}



\ECSwitch


\ECHead{Electronic Companion}


\section{Proofs of Solution Quality Bounds}


\subsection{Proof of Theorem \ref{thm:non-parametric_bound}}

\proof{Proof.}
\underline{We start by defining notations}. For any leader decision $\bfx\in \mX$, let
\begin{equation*}
    \hat{F}_\mT^{k\textrm{NN}}(\bfx) := 
    f(\bfx) 
        + \sum_{t\in \mT} q^{t} G^{t}(\bfx)
        + \sum_{s\in \mS\backslash\mT}
            \sum_{t\in \mT_k(\bff^s)} 
                \frac{q^{s}}{k}
                G^t(\bfx)
\end{equation*}
denote the objective function of the $k$NN-augmented model~\eqref{prob:non-parametric_approx}. Recall that $F(\bfx)$ is the objective function of the original bilevel model~\eqref{prob:general} that considers the full follower set $\mS$. For any $\bfx \in \mX$ and $s\in \mS$, let $\eta^s(\bfx) := \bbE\left[G^s(\bfx) \,|\, \bff^s \right]$ be the expected cost of follower $s$. Let $\bar{Q} := \max_{s\in \mS} q^s$ denote the maximal follower weight.

\underline{We first decompose the approximation error of $\hat{F}^{k\textrm{NN}}_\mT(\bfx)$ to $F(\bfx)$.} For any $\bfx\in \mX$, 
\begin{align*}
& F(\bfx) - \hat{F}_{\mT}^{k\textrm{NN}}(\bfx) \\
=
& \sum_{s\in \mS\backslash\mT}
    q^s
    \left(
        G^s(\bfx)
        - \sum_{t\in \mT_k(\bff^s)}
            \frac{1}{k}
            G^t(\bfx)
    \right) \\
=  
&\sum_{s\in \mS\backslash\mT}
    q^s
    \left(
        G^s(\bfx) 
        + \eta^s(\bfx) 
        - \eta^s(\bfx)
        - \sum_{t\in \mT_k(\bff^s)}
            \frac{1}{k}
            G^t(\bfx)
        + \sum_{t\in \mT_k(\bff^s)}
            \frac{1}{k}
            \eta^t(\bfx)
        - \sum_{t\in \mT_k(\bff^s)}
            \frac{1}{k}
            \eta^t(\bfx)
    \right) \\
= 
&
\sum_{s\in \mS}
\sum_{t\in \mT_k(\bff^s)}
    \frac{q^s}{k}
    \left[
        \eta^s(\bfx) - \eta^t(\bfx)
    \right]
+ \sum_{s\in \mS\backslash\mT}
    q^s 
    \left[
        G^s(\bfx) - \eta^s(\bfx)
    \right]
+ \sum_{t\in \mT}
    \sum_{s\in \mS_k(t)}
    \frac{q^s}{k}
    \left[
        \eta^t(\bfx) - G^t(\bfx)
    \right],
& 
\end{align*}
where $\mS_k(t) \subseteq \mS$ represents the set of followers in $\mS\backslash\mT$ whose $k$-nearest neighbors in set $\mT$ include follower $t$. 

\underline{Next, we relate the approximation error to decision quality.} We have
\begin{align*}
F(\bfx^{k\textrm{NN}}_\mT) 
    - F(\bfx^*) 
= 
& F(\bfx^{k\textrm{NN}}_\mT) 
    - \hat{F}_{\mT}^{k\textrm{NN}}(\bfx^{k\textrm{NN}}_\mT)
    + \hat{F}_{\mT}^{k\textrm{NN}}(\bfx^{k\textrm{NN}}_\mT)
    - F(\bfx^*)
    + \hat{F}_{\mT}^{k\textrm{NN}}(\bfx^*)
    - \hat{F}_{\mT}^{k\textrm{NN}}(\bfx^*) \\
\leq
& \left[ 
    F(\bfx^{k\textrm{NN}}_\mT) 
    - \hat{F}_{\mT}^{k\textrm{NN}}(\bfx^{k\textrm{NN}}_\mT)
\right]
+
\left[ 
    \hat{F}_{\mT}^{k\textrm{NN}}(\bfx^*)
    - F(\bfx^*)
\right]
\end{align*} 

Based on the approximation error decomposition above, we have
\begin{align*}
& F(\bfx^{k\textrm{NN}}_\mT) - F(\bfx^*) \\
\leq 
& 
\underbrace{\sum_{s\in \mS}
\sum_{t\in \mT_k(\bff^s)}
    \frac{q^s}{k}
    \left[
        \eta^s(\bfx^{k\textrm{NN}}_\mT) - \eta^t(\bfx^{k\textrm{NN}}_\mT)
    \right]}_\textrm{(1)}
+ \underbrace{\sum_{s\in \mS\backslash\mT}
    q^s
    \left[
        G^s(\bfx^{k\textrm{NN}}_\mT) 
        - \eta^s(\bfx^{k\textrm{NN}}_\mT)
    \right]}_\textrm{(2)} \\
& 
+ \underbrace{\sum_{t\in \mT}
   \sum_{s\in \mS_k(t)}
    \frac{q^s}{k}
    \left[
        \eta^t(\bfx^{k\textrm{NN}}_\mT)
        - G^t(\bfx^{k\textrm{NN}}_\mT) 
    \right]}_\textrm{(3)} 
+ \underbrace{\sum_{s\in \mS}
\sum_{t\in \mT_k(\bff^s)}
    \frac{q^s}{k}
    \left[
        \eta^t(\bfx^*) - \eta^s(\bfx^*)
    \right]}_\textrm{(4)} \\
& 
+ \underbrace{\sum_{s\in \mS\backslash\mT}
    q^s
    \left[ 
        \eta^s(\bfx^*) - G^s(\bfx^*) 
    \right]}_\textrm{(5)}
+ \underbrace{\sum_{t\in \mT}
    \sum_{s\in \mS_k(t)}
    \frac{q^s}{k}
    \left[
        G^t(\bfx^*) 
        - \eta^t(\bfx^*)
    \right]}_\textrm{(6)}
\end{align*}

According to Assumption \ref{asm:lipschitz_cost}, we have
\begin{align*}
(1) + (4)
\leq 
&
\sum_{s\in \mS}
\sum_{t\in \mT_k(\bff^s)}
    \frac{\bar{Q}}{k}
    \left|
        \eta^s(\bfx^{k\textrm{NN}}_\mT) - \eta^t(\bfx^{k\textrm{NN}}_\mT)
    \right|
+ 
\sum_{s\in \mS}
\sum_{t\in \mT_k(\bff^s)}
    \frac{\bar{Q}}{k}
    \left|
        \eta^s(\bfx^*) - \eta^t(\bfx^*)
    \right|\\
\leq
&
\sum_{s\in \mS}
\sum_{t\in \mT_k(\bff^s)}
    \frac{2\mu\bar{Q}}{k}
    d(\bff^s, \bff^t).
\end{align*}

Regarding (2), (3), (5), and (6), for any $s\in \mS\backslash\mT$, we can treat $q^s \left[G^s(\bfx^{k\textrm{NN}}_\mT) - \eta^s(\bfx^{k\textrm{NN}}_\mT)\right]$ and $q^s \left[G^s(\bfx^*) - \eta^s(\bfx^*)\right]$ as random variables with mean zero and are bounded in an interval of width $2\bar{Q}\bar{G}$. Similarly, for any $t\in \mT$, we can treat $\sum_{s\in \mS_k(t)} q^s \left[ G^t(\bfx) - \eta^t(\bfx) \right] / k$ and $\sum_{s\in \mS_k(t)} q^s \left[\eta^t(\bfx) - G^t(\bfx) \right] / k$ as random variables with mean zero and that are bounded in an interval of width $2\bar{Q} \bar{G} m^t_{k, \mS\backslash\mT} / k$. Based on Assumption \ref{asm:iid} and the Hoeffding inequality \citep{hoeffding1994probability}, we have, with probability at least $1 - \gamma$,
\begin{equation*}
    (2) + (3) + (5) + (6)
    \leq 
    \sqrt{
        4\bar{Q}^2\bar{G}^2
        \log\left(\frac{1}{\gamma}\right)
        \left(
            |\mS\backslash\mT|
            + \sum_{t\in \mT}
                \left(\frac{m^t_{k, \mS\backslash \mT}}{k}\right)^2
        \right)
    },
\end{equation*}
where $m^t_{k, \mS\backslash \mT}$ is the number of followers in $\mS\backslash\mT$ whose $k$-nearest neighbors in $\mT$ include follower $t\in \mT$.

Combining the bounds on (1) + (4) and (2) + (3) + (5) + (6), we have, with probability at least $1 - \gamma$,
\begin{equation*}
F(\bfx^{k\textrm{NN}}_\mT) - F(\bfx^*)
\leq 
\sum_{s\in \mS}
\sum_{t\in \mT_k(\bff^s)}
    \frac{2\mu\bar{Q}}{k}
    d(\bff^s, \bff^t)
+
\sqrt{
    4\bar{Q}^2\bar{G}^2
    \log\left(\frac{1}{\gamma}\right)
    \left(
        |\mS\backslash\mT|
        + \sum_{t\in \mT}
            \left(\frac{m^t_{k, \mS\backslash \mT}}{k}\right)^2
    \right)
}.
\end{equation*}
\Halmos
\endproof

\subsection{Proof of Theorem \ref{thm:reg_bound}}

\underline{We start by defining notations}. Given a parametric regression model $P(\cdot; \bstheta)$, for any leader decision $\bfx\in \mX$, let
\begin{align*}
    \hat{F}_\mT^{\textrm{PR}}(\bfx) := 
    \min_{\bstheta\in \bsTheta} \ 
    & f(\bfx) 
        + \sum_{t\in \mT} q^{t} G^{t}(\bfx)
        + \sum_{s\in \mS\backslash\mT}
            P(\bff^s; \bstheta) \\
    \textrm{s.t.} \quad
    & \sum_{t\in \mT} 
            m^t_{1, \mS\backslash\mT}
            \left| G^t(\bfx) - P(\bff^t; \bstheta) \right| 
        \leq \bar{L}
\end{align*}
denote the objective function of the parametric regression-augmented model~\eqref{prob:reg_approx}. Recall that $F(\bfx)$ is the objective function of the original bilevel model~\eqref{prob:general} that considers the full follower set $\mS$. For any $\bfx \in \mX$ and $s\in \mS$, let $\eta^s(\bfx) := \bbE\left[G^s(\bfx) \,|\, \bff^s \right]$ be the expected cost of follower $s$. Let $\bar{Q} := \max_{s\in \mS} q^s$ denote the maximal follower weight. We denote by $\bstheta^\textrm{PR}_{\mT}$ and $\bstheta^*_{\mT}$ the optimal solutions to the optimization models solved for evaluating $\hat{F}_\mT^{\textrm{PR}}(\bfx^\textrm{PR})$ and $\hat{F}_\mT^{\textrm{PR}}(\bfx^*)$, respectively. Let $\nu(\bff^s)$ be the nearest neighbor of $\bff^s$ in $\mT$ for any $s\in \mS\backslash\mT$. 

\underline{We first decompose the approximation error of $\hat{F}^{\textrm{PR}}_\mT(\bfx)$ to $F(\bfx)$.} For any $\bfx\in \mX$, we have
\begin{align*}
& F(\bfx) - \hat{F}_{\mT}^{\textrm{PR}}(\bfx) \\
=
& \sum_{s\in \mS\backslash\mT}
    q^s
    \left[
        G^s(\bfx)
        - P(\bff^s; \bstheta)
    \right] \\
=
& 
\sum_{s\in \mS\backslash\mT}
    q^s
    \left[
        G^s(\bfx)
        - \eta^s(\bfx)
    \right] 
+
\sum_{s\in \mS\backslash\mT}
    q^s
    \left[
        \eta^{s}(\bfx)
        - \eta^{\nu(\bff^s)}(\bfx)
    \right] 
+
\sum_{s\in \mS\backslash\mT}
    q^s
    \left[
        \eta^{\nu(\bff^s)}(\bfx)
        - G^{\nu(\bff^s)}(\bfx)
    \right] \\
&
+
\sum_{s\in \mS\backslash\mT}
    q^s
    \left[
        G^{\nu(\bff^s)}(\bfx)
        - P(\bff^{\nu(\bff^s)}; \bstheta)
    \right]
+
\sum_{s\in \mS\backslash\mT}
    q^s
    \left[
        P(\bff^{\nu(\bff^s)}; \bstheta)
        - P(\bff^{s}; \bstheta)
    \right] \\ 
=
& 
\sum_{s\in \mS\backslash\mT}
    q^s
    \left[
        G^s(\bfx)
        - \eta^s(\bfx)
    \right] 
+
\sum_{s\in \mS\backslash\mT}
    q^s
    \left[
        \eta^{s}(\bfx)
        - \eta^{\nu(\bff^s)}(\bfx)
    \right] 
+
\sum_{t \in \mT}
\sum_{s\in \mS_1(t)}
    q^s
    \left[
        \eta^{t}(\bfx)
        - G^{t}(\bfx)
    \right] \\
&
+
\sum_{t\in \mT}
    \sum_{s\in \mS_1(t)}
    q^s
    \left[
        G^t(\bfx)
        - P(\bff^t; \bstheta)
    \right]
+
\sum_{s\in \mS\backslash\mT}
    q^s
    \left[
        P(\bff^{\nu(\bff^s)}; \bstheta)
        - P(\bff^{s}; \bstheta)
    \right]
\end{align*}
where $\mS_k(t) \subseteq \mS$ represents the set of followers in $\mS\backslash\mT$ whose $k$-nearest neighbors in set $\mT$ include follower $t$. 

\underline{Next, we relate the approximation error to decision quality.} We have
\begin{align*}
& F(\bfx^{\textrm{PR}}_\mT) 
    - F(\bfx^*) \\
= 
& F(\bfx^{\textrm{PR}}_\mT) 
    - \hat{F}_{\mT}^{\textrm{PR}}(\bfx^{\textrm{PR}}_\mT)
    + \hat{F}_{\mT}^{\textrm{PR}}(\bfx^{\textrm{PR}}_\mT)
    - F(\bfx^*)
    + \hat{F}_{\mT}^{\textrm{PR}}(\bfx^*)
    - \hat{F}_{\mT}^{\textrm{PR}}(\bfx^*) \\
\leq
& \left[ 
    F(\bfx^{\textrm{PR}}_\mT) 
    - \hat{F}_{\mT}^{\textrm{PR}}(\bfx^{\textrm{PR}}_\mT)
\right]
+
\left[ 
    \hat{F}_{\mT}^{\textrm{PR}}(\bfx^*)
    - F(\bfx^*)
\right]
\end{align*}
The inequality holds because $\bfx^{\textrm{PR}}_\mT$ is the optimal solution to model~\eqref{prob:reg_approx}. 

Based on the approximation error decomposition above, we have
\begin{align*}
& F(\bfx^{\textrm{PR}}_\mT) - F(\bfx^*) \\
\leq 
&
\underbrace{\sum_{s\in \mS\backslash\mT}
    q^s
    \left[
        G^s(\bfx^{\textrm{PR}}_{\mT})
        - \eta^s(\bfx^{\textrm{PR}}_{\mT})
    \right]}_\textrm{(1)}
+
\underbrace{\sum_{s\in \mS\backslash\mT}
    q^s
    \left[
        \eta^{s}(\bfx^{\textrm{PR}}_{\mT})
        - \eta^{\nu(\bff^s)}(\bfx^{\textrm{PR}}_{\mT})
    \right]}_\textrm{(2)}
+
\underbrace{\sum_{t \in \mT}
\sum_{s\in \mS_1(t)}
    q^s
    \left[
        \eta^{t}(\bfx^{\textrm{PR}}_{\mT})
        - G^{t}(\bfx^{\textrm{PR}}_{\mT})
    \right]}_\textrm{(3)} \\
&
+
\underbrace{\sum_{t\in \mT}
    \sum_{s\in \mS_1(t)}
    q^s
    \left[
        G^t(\bfx^{\textrm{PR}}_{\mT})
        - P(\bff^t; \bstheta^{\textrm{PR}}_{\mT})
    \right]}_\textrm{(4)}
+
\underbrace{\sum_{s\in \mS\backslash\mT}
    q^s
    \left[
        P(\bff^{\nu(\bff^s)}; \bstheta^{\textrm{PR}}_{\mT})
        - P(\bff^{s}; \bstheta^{\textrm{PR}}_{\mT})
    \right]}_\textrm{(5)} \\
+ 
& 
\underbrace{\sum_{s\in \mS\backslash\mT}
    q^s
    \left[
        G^s(\bfx^*)
        - \eta^s(\bfx^*)
    \right]}_\textrm{(6)} 
+
\underbrace{\sum_{s\in \mS\backslash\mT}
    q^s
    \left[
        \eta^{s}(\bfx^*)
        - \eta^{\nu(\bff^s)}(\bfx^*)
    \right]}_\textrm{(7)} 
+
\underbrace{\sum_{t \in \mT}
\sum_{s\in \mS_1(t)}
    q^s
    \left[
        \eta^{t}(\bfx^*)
        - G^{t}(\bfx^*)
    \right]}_\textrm{(8)} \\
&
+
\underbrace{\sum_{t\in \mT}
    \sum_{s\in \mS_1(t)}
    q^s
    \left[
        G^t(\bfx^*)
        - P(\bff^t; \bstheta^*_{\mT})
    \right]}_\textrm{(9)}
+
\underbrace{\sum_{s\in \mS\backslash\mT}
    q^s
    \left[
        P(\bff^{\nu(\bff^s)}; \bstheta^*_{\mT})
        - P(\bff^{s}; \bstheta^*_{\mT})
    \right]}_\textrm{(10)}
\end{align*}

According to Assumption \ref{asm:lipschitz_cost}, we have
\begin{align*}
(2) + (7)
\leq 
\sum_{s\in \mS}
    \bar{Q}
    \left|
        \eta^{s}(\bfx^{\textrm{PR}}_{\mT})
        - \eta^{\nu(\bff^s)}(\bfx^{\textrm{PR}}_{\mT})
    \right|
+ 
\sum_{s\in \mS}
    \bar{Q}
    \left|
        \eta^{s}(\bfx^*)
        - \eta^{\nu(\bff^s)}(\bfx^*)
    \right|
\leq
\sum_{s\in \mS}
    2\mu\bar{Q}
    d(\bff^s, \bff^t).
\end{align*}

According to Assumption \ref{asm:lipschitz_pred}, we have
\begin{align*}
(5) + (10) 
\leq 
& \sum_{s\in \mS\backslash\mT}
    \bar{Q}
    \left|
        P(\bff^{\nu(\bff^s)}; \bstheta^*_{\mT})
        - P(\bff^{s}; \bstheta^*_{\mT})
    \right|
    +
    \sum_{s\in \mS\backslash\mT}
    \bar{Q}
    \left|
        P(\bff^{\nu(\bff^s)}; \bstheta^*_{\mT})
        - P(\bff^{s}; \bstheta^*_{\mT})
    \right| \\
\leq 
& \sum_{s\in \mS\backslash\mT}
    2\lambda \bar{Q} 
    d(\bff^s, \bff^t).
\end{align*}

Since both $(\bfx^*, \bstheta_\mT^*)$ and $(\bfx^\textrm{PR}_\mT, \bstheta_\mT^\textrm{PR})$ are both feasible solutions to model~\eqref{prob:reg_approx}, we have
\begin{align*}
(4) + (9) 
\leq 
& \sum_{t \in \mT}
    \bar{Q}
    m_{1, \mS\backslash\mT}^t
    \left|
        G^t(\bfx^{\textrm{PR}}_{\mT})
        - P(\bff^t; \bstheta^{\textrm{PR}}_{\mT})
    \right|
+ \sum_{t\in \mT}
    \bar{Q}
    m_{1, \mS\backslash\mT}^t
    \left|
        G^t(\bfx^*)
        - P(\bff^t; \bstheta^*_{\mT})
    \right| \\
\leq 
& 2\bar{Q} \bar{L} \mathds{1}(\mT \subset \mS)
\end{align*}

The last inequality holds because we need the training loss constraint only when $\mT \subset \mS$.

Regarding (1), (3), (6), and (8), for any $s\in \mS\backslash\mT$, we can treat $q^s \left[G^s(\bfx^{\textrm{PR}}_\mT) - \eta^s(\bfx^{\textrm{PR}}_\mT)\right]$ and $q^s \left[G^s(\bfx^*) - \eta^s(\bfx^*)\right]$ as random variables with mean zero and are bounded in an interval of width $2\bar{Q}\bar{G}$. Similarly, for any $t\in \mT$, we can treat $\sum_{s\in \mS_1(t)}q^s \left[ G^t(\bfx) - \eta^t(\bfx) \right]$ and $\sum_{s\in \mS_1(t)} q^s \left[\eta^t(\bfx) - G^t(\bfx) \right]$ as random variables with mean zero and that are bounded in an interval of width $2 \bar{Q} \bar{G} m^t_{1, \mS\backslash\mK}$. Based on Assumption \ref{asm:iid} and the Hoeffding inequality \citep{hoeffding1994probability}, we have, with probability at least $1 - \gamma$,
\begin{equation*}
    (1) + (3) + (6) + (8)
    \leq 
    \sqrt{
        4\bar{Q}^2\bar{G}^2
        \log\left(\frac{1}{\gamma}\right)
        \left(
            |\mS\backslash\mT|
            + \sum_{t\in \mT}
                \left(m^t_{1, \mS\backslash \mT}\right)^2
        \right)
    },
\end{equation*}
where $m^t_{k, \mS\backslash \mT}$ is the number of followers in $\mS\backslash\mT$ whose $k$-nearest neighbors in $\mT$ include follower $t\in \mT$.

Combing the bounds on (2) + (7), (5) + (10), (4) + (9), and (1) + (3) + (6) + (8), we have with probability at least $1 - \gamma$,
\begin{align*}
    & F(\bfx^{PR}_\mT) - F(\bfx^*) \\
    \leq 
    & 2 \bar{Q}\bar{L} \mathds{1}(\mT \subset \mS)
    + 2 \bar{Q} (\lambda + \mu)
    \sum_{s\in \mS\backslash\mT} 
        d_{\mF}(\bff^s, \bff^{\nu(s)})
    + \sqrt{4\bar{Q}^2\bar{G}^2\left[|\mS\backslash\mT| + \sum_{t\in\mT}(m^t_{1, \mS\backslash\mT})^2 \right]\log(1/\gamma)}.
\end{align*}
\Halmos
\endproof

\section{Proofs of Bound Tightness Results}
\subsection{Proof of Theorem \ref{thm:knn_tightness}}

\underline{We first define some notations that will be used in this proof.} We rewrite the bound in Theorem \ref{thm:non-parametric_bound} as follows.
\begin{equation*}
    E^{1\textrm{NN}}_m(\mT) = 
    2 \mu \bar{Q} 
    B^{1\textrm{NN}}(\mT)
    + 2\bar{Q} \bar{G} V^{1\textrm{NN}}(\mT),
\end{equation*}
where
\begin{equation*}
    B^{1\textrm{NN}}(\mT) 
    := \sum_{s\in \mS}
        \sum_{t\in \mT_k(\bff^s)}
            d(\bff^s, \bff^t)
\end{equation*}
and
\begin{equation*}
    V^{1\textrm{NN}}(\mT)
    := \sqrt{\left[|\mS\backslash\mT|
    + \sum_{t\in\mT}(m^t_{1, \mS\backslash\mT})^2\right]\log(1/\gamma)}
\end{equation*}
indicate the bias and variance terms, respectively. Next, we derive the limits of $B^{1\textrm{NN}}(\mT_{p, d}) / m$ and $V^{1\textrm{NN}}(\mT_{p, d}) / m$, separately.

\underline{Limit of $B^{1\textrm{NN}}(\mT_{p, d}) / m$}. Observing that $\mT_{p, d}$ is the optimal solution to a balanced $p$-median problem, according to Theorem 1.1 in \citet{mcgivney1999asymptotics}, we have
\begin{equation*}
    \lim_{m\rightarrow \infty} 
        \frac{1}{m^{(\xi-1) / \xi}}
        B^{1\textrm{NN}}(\mT_{p, d}) 
    = C_{d, \xi} \int_{[0, 1]^\xi} \sigma^{(\xi-1)/\xi}(\bff) df.
\end{equation*}
where $C_{d, \xi} > 0$ is a constant that depends on $d$ and $\xi$. As estimated in Theorem 12 from \citet{carlsson2022continuous}, $C_{d, \xi}$ satisfies
\begin{equation*}
    C_{d, \xi} \leq \frac{2}{3}\sqrt{d}.
\end{equation*}
Since $d\leq \beta \lceil m/p \rceil$ and $p = \max \{1, \alpha m^{(\xi -1) / \xi} \}$, we have
\begin{equation*}
    \lim_{m\rightarrow \infty}
        \frac{1}{m}
        B^{1\textrm{NN}}(\mT_{p, d})
        = 0.
    \label{eq:lim_knn_first}
\end{equation*}

\underline{Limit of $V^{1\textrm{NN}}(\mT_{p, d}) / m$}. It is easy to verify that 
\begin{equation*}
    0 
    \leq 
    V^{1\textrm{NN}}(\mT_{p, d})
    \leq 
    \bar{V}^{1\textrm{NN}}(\mT_{p, d}),
\end{equation*}
where
\begin{equation*}
    \bar{V}^{1\textrm{NN}}(\mT)
    =
    \sqrt{
        \left[
            |\mS|
            - |\mT|
            + |\mT| d^2
        \right]
        \log(1/\gamma)
    }.
\end{equation*}
The first inequality is trivial, while the second inequality holds because $m_{1, \backslash\mT}^t = |\mS(t)| \leq d$ for all $t\in \mT$ which come from the constraints in Problem~\eqref{prob:follower_select_knn}. Therefore, according to the Squeeze Theorem, to show 
\begin{equation*}
    \lim_{m\rightarrow \infty} 
        \frac{1}{m} 
        V^{1\textrm{NN}}(\mT_{p, d}) = 0,
\end{equation*}
it suffices to show that
\begin{equation*}
    \lim_{m\rightarrow \infty} 
        \frac{1}{m} 
        \bar{V}^{1\textrm{NN}}(\mT_{p, d}) = 0.
\end{equation*}

We have 
\begin{align*}
0
\leq 
\frac{1}{m} 
    \bar{V}^{1\textrm{NN}}(\mT_{p, d})
\leq
\sqrt{\left[\frac{1}{|\mS|} - \frac{\alpha}{|\mS|^{(\xi + 1)/\xi}} + \frac{\beta^2}{\alpha |\mS|^{(\xi -1) / \xi}} \right]
        \log(1/\gamma)}
\end{align*}

According to the Squeeze theorem, we have
\begin{equation*}
     \lim_{m\rightarrow \infty} 
        \frac{1}{m} 
        V^{1\textrm{NN}}(\mT_{p, d}) = 0
\end{equation*}

\underline{Combining the limits of $V^{1\textrm{NN}}(\mT_{p, d}) / m$ and $B^{1\textrm{NN}}(\mT_{p, d}) / m$}, we have
\begin{equation*}
    \lim_{m\rightarrow \infty} 
        \frac{1}{m} 
        E^{1\textrm{NN}}_m(\mT_{p, d}) = 0.
\end{equation*}
\Halmos
\endproof

\subsection{Proof of Theorem \ref{thm:reg_tightness}}

\underline{We first define some notations that will be used in this proof.} We rewrite the bound in Theorem \ref{thm:reg_bound} as follows.
\begin{equation*}
    E^{\textrm{PR}}_m(\mT) = 
    2\bar{Q}\bar{L} \mathds{1}(\mT \subset \mS)
    + 2 \bar{Q} (\mu + \lambda)
        B^{1\textrm{NN}}(\mT)
    + 2\bar{Q} \bar{G} V^{1\textrm{NN}}(\mT),
\end{equation*}
where
\begin{equation*}
    B^{1\textrm{NN}}(\mT) 
    := \sum_{s\in \mS}
        \sum_{t\in \mT_k(\bff^s)}
            d(\bff^s, \bff^t)
\end{equation*}
and
\begin{equation*}
    V^{1\textrm{NN}}(\mT)
    := \sqrt{\left[|\mS\backslash\mT|
    + \sum_{t\in\mT}(m^t_{1, \mS\backslash\mT})^2\right]\log(1/\gamma)}
\end{equation*}
indicate the bias and variance terms, respectively.

As shown in the proof of Theorem \ref{thm:knn_tightness}, we know that 
\begin{equation*}
    \lim_{m\rightarrow \infty}
        \frac{1}{m}
        B^{1\textrm{NN}}(\mT) 
    = 0,
\end{equation*}
and 
\begin{equation*}
    \lim_{m\rightarrow \infty}
        \frac{1}{m}
        V^{1\textrm{NN}}(\mT) 
    = 0.
\end{equation*}

Since the first term $2\bar{Q}\bar{L} \mathds{1}(\mT \subset \mS)$ is a fixed finite constant, we have 
\begin{equation*}
    \lim_{m\rightarrow \infty}
        \frac{1}{m}
        E^{\textrm{PR}}_m(\mT) 
    = 0.
\end{equation*}
\Halmos
\endproof

\section{Formulation of MaxANDP} \label{app:maxandp_formulation}
In this section, we present the full formulation of MaxANDP using a piecewise linear impedance function $g$ and shortest-path routing problems as introduced in Section \ref{subsec:prob_formulation}.


Since $g$ is a decreasing function of travel time and the objectives of the followers are to minimize travel time, the objectives of the leader and followers are aligned. Therefore, the optimality conditions~\eqref{cyc_general:opt_route} can be replaced with the routing constraints~\eqref{od_routing:balance}--\eqref{od_routing:node_design}, resulting in a single-level formulation:
\begin{subequations}
\label{prob:maxandp_piece_single}
\begin{align}
    \underset{\bfx, \bfy, \bfz}{\textrm{maximize}} \quad 
        &\sum_{(o, d) \in \mS}
        q^{od} g(\bfy^{od}) \\
    \textrm{subject to} \quad
    \label{full:budget_edge}
    & \bfc^\intercal \bfx \leq B_\text{edge} \\
    \label{full:budget_node}
    & \bfb^\intercal \bfz \leq B_\text{node} \\
    \label{full:flow_balance}
    & \textbf{A}\bfy^{od} = \bfe^{od},
        \quad \forall (o, d) \in \mS\\
    \label{full:edge_design}
    & y^{od}_{ij} \leq x_{ij},
        \quad \forall (i,j)\in \mE_{h}, \ (o, d)\in \mS\\
    \label{full:node_design}
    & y^{od}_{ij} \leq x_{wl} + z_i, 
        \quad \forall i \in \mN_h, 
        \ (i,j)\in \mE^-_{h}(i), 
        \ (w,l)\in \mE_{h}(i),
        \ (o, d)\in \mS\\
    \label{full:domain_y}
    & 0 \leq y^{od}_{ij} \leq 1, 
        \quad \forall (i,j)\in \mE
        \ (o, d) \in \mS \\
    \label{full:domain_x}
    & \bfx\in \{0, 1\}^{|\mE_h|} \\
    \label{full:domain_z}
    & \bfz\in \{0, 1\}^{|\mN_h|}.
\end{align}
\end{subequations}
It is known that the parameter matrix of the flow balance constraints~\eqref{full:flow_balance} is totally unimodular. Moreover, for any fixed feasible $\bfx$ and $\bfz$, each row in the parameter matrix of constraints \eqref{full:edge_design}--\eqref{full:domain_y} has one ``1'' with all other entries being zero. Therefore, for any fixed $\bfx$ and $\bfz$, the parameter matrix of constraints~\eqref{full:flow_balance}--\eqref{full:domain_y} is totally unimodular. We thus can discard the integrality constraints on $\bfy^{od}$ and treat them as continuous decision variables bounded in $[0, 1]$. The only thing left is to linearize function $g$, which depends on the values of $\beta_1$ and $\beta_2$.

\textbf{Concave impedance function.} When $\beta_1 \leq \beta_2$, $g$ can be treated as a concave function of travel time $\bft^\intercal \bfy^{od}$ on $[0, T_2]$. We introduce continuous decision variables $v^{od} \in \bbR_+$, for any OD pairs $(o, d)\in \mS$, representing the accessibility of OD pairs $(o, d)$. Problem~\eqref{prob:maxandp_piece_single} can then be written as
\begin{subequations}
\begin{align}
    \underset{\bfv, \bfx, \bfy, \bfz}{\textrm{maximize}} \quad 
        &\sum_{(o, d) \in \mS}
        q^{od}v^{od} \\
    \notag \textrm{subject to} \quad
    & \eqref{full:budget_edge}\text{--}\eqref{full:node_design} \\
    & v^{od} \leq \alpha_1 - \beta_1 \bft^\intercal \bfy^{od},
        \quad \forall (o, d) \in \mS \\
    & v^{od} \leq \alpha_2 - \beta_2 \bft^\intercal \bfy^{od},
        \quad \forall (o, d) \in \mS \\
    & v^{od} \geq 0, \quad \forall (o, d) \in \mS \\
    \notag
    & \eqref{full:domain_y}\text{--}\eqref{full:domain_z}
\end{align}
\end{subequations}
where $\alpha_1=1$ and $\alpha_2 = 1 + (\beta_2 - \beta_1) T_1$ are the intercepts of the linear functions in $[0, T_1)$ and $[T_1, T_2)$ respectively.

\textbf{Convex impedance function}. When $\beta_1 > \beta_2$, $g$ can be treated as a strictly convex function of travel time $\bft^\intercal \bfy^{od}$ on $[0, T_2]$. We introduce binary decision variables $r^{od}$, for any OD pairs $(o, d)\in \mS$, representing if the travel time of pair $(o, d)$ is in $[0, T_1]$ ($=1$) or not ($=0$). We introduce continuous decision variabel $u^{od}_t \in \bbR_+$, for any OD pairs $(o, d)\in \mS$ and any of the two travel time intervals $t\in \{1, 2\}$. Problem~\eqref{prob:maxandp_piece_single} can then be written as 
\begin{subequations}
\begin{align}
    \underset{\bfr, \bfu, \bfx, \bfy, \bfz}{\textrm{maximize}} \quad 
        &\sum_{(o, d) \in \mS}
        \sum_{t=1}^{2}
            q^{od} (\alpha_t - \beta_t u^{od}_t) \\
    \notag \textrm{subject to} \quad
    & \eqref{full:budget_edge}\text{--}\eqref{full:node_design} \\
    & u^{od}_1 \leq T_1 r^{od}, 
        \quad \forall (o, d) \in \mS \\
    & T_1 (1 - r^{od})\leq u^{od}_2 \leq T (1 - r^{od}), 
        \quad \forall (o, d) \in \mS \\
    & u^{od}_1 + u^{od}_2 = \bft^\intercal \bfy^{od}, 
        \quad \forall (o, d) \in \mS \\
    & u^{od}_t \geq 0, 
        \quad \forall (o, d) \in \mS, \ t\in \{1, 2\} \\
    & r^{od} \in \{0, 1\}, 
        \quad \forall (o, d) \in \mS \\
    & \eqref{full:domain_y}\text{--}\eqref{full:domain_z}.
\end{align}
\end{subequations}

\section{Solution Method for MaxANDP} \label{app:benders}

In this section, we present a Benders decomposition algorithm that takes advantage of the block structure of Problem~\eqref{prob:maxandp_piece_single} (i.e. the routing problems of the OD pairs are independent of each other) along with its acceleration strategies. This algorithm is applicable to both convex and concave impedance functions, and can be easily extended to solving the ML-augmented model of MaxANDP. For ease of presentation, we treat $g$ as a function for travel time $\bft^\intercal \bfy^{od}$ in this section. 

\subsection{Benders Decomposition}
We start by introducing the Benders reformulation of Problem~\eqref{prob:maxandp_piece_single}. We introduce continuous decision variables $\tau^{od}$, for any OD pair $(o, d)\in \mS$, indicating the travel time from $o$ to $d$ using the low-stress network. We then re-written Problem~\eqref{prob:maxandp_piece_single} as
\begin{subequations}
\label{prob:maxan_benders_reformulation}
\begin{align}
    \underset{\bfx, \bfy, \bfz, \boldsymbol{\tau}}{\textrm{maximize}} \quad 
        &\sum_{(o, d) \in \mS}
        q^{d} g(\tau^{od}) \\
    \textrm{subject to} \quad
    \notag
    & \eqref{full:budget_edge}, \eqref{full:budget_node}, \eqref{full:domain_x}, \eqref{full:domain_z} \\
    \label{reform:time_primal}
    & \tau^{od} \geq 
        \min_{\bfy^{od}}
            \left\{ 
                \bft^\intercal \bfy^{od}
            \,\middle|\,
                \eqref{full:flow_balance}-\eqref{full:domain_y}
            \right\}, 
        \quad \forall (o, d) \in \mS.
\end{align}
\end{subequations}
For each $(o, d)\in \mS$, We associate unbounded dual variables $\boldsymbol{\lambda}^{od}$ with constraints \eqref{full:flow_balance} and non-negative dual variables $\boldsymbol{\theta}^{od}$, $\boldsymbol{\delta}^{od}$, and $\boldsymbol{\pi}^{od}$ with constraints \eqref{full:edge_design}--\eqref{full:domain_y}, respectively. Given any network design $(\bfx, \bfz)$, we formulate the dual of the routing problem as
\begin{subequations}
\label{prob:benders_sub_dual}
\begin{align}
    \underset{ \boldsymbol{\theta}, \boldsymbol{\delta}, \boldsymbol{\pi} 
    \geq 0,\ \boldsymbol{\lambda}}{\textrm{maximize}} \quad 
    & -\lambda_{d}^{od} + \lambda_{o}^{od}
        - \sum_{(i,j)\in \mE_h} 
            x_{ij}\theta_{ij}^{od}
        - \sum_{i\in \mN_h}
            \sum_{(i,j)\in \mE^-_{h}(i)}
            \sum_{(w,l)\in \mE_{h}(i)}
            (x_{wl} + z_i) \delta_{ijwl}^{od} 
        - \sum_{(i, j)\in \mE} \pi^{od}_{ij} \\
    \textrm{subject to} \quad
    & - \lambda_{j}^{od} 
        + \lambda_{i}^{od}
        - \mathds{1}\left[(i, j) \in \mE_h \right] \theta_{ij}^{od}
        - \mathds{1}(i\in \mN_h) 
            \sum_{(w, l)\in \mE_{h}(i)}
            \delta_{ijwl}^{od} 
        - \pi^{od}_{ij}
        \leq t_{ij},
        \quad \forall (i, j) \in \mE.
\end{align}
\end{subequations}
Since the routing problem associated with each OD pair is always feasible and bounded, its dual~\eqref{prob:benders_sub_dual} is also feasible and bounded. Let $\Pi^{od}$ denote the set of extreme points of problem \eqref{prob:benders_sub_dual}. According to the duality theory, constraints \eqref{reform:time_primal} can be replace by
\begin{multline}
\label{con:cut_plane}
    \tau^{od} \geq 
        -\lambda_{d}^{od} + \lambda_{o}^{od}
        - \sum_{(i,j)\in \mE_h} 
            x_{ij}\theta_{ij}^{od}
        - \sum_{i\in \mN_h}
            \sum_{(i,j)\in \mE^-_{h}(i)}
            \sum_{(w,l)\in \mE_{h}(i)}
            (x_{wl} + z_i) \delta_{ijwl}^{od} 
        - \sum_{(i, j)\in \mE} \pi^{od}_{ij} \\
    \quad \forall
        \boldsymbol{\lambda}^{od}, \boldsymbol{\theta}^{od}, \boldsymbol{\delta}^{od}, \boldsymbol{\pi}^{od} \in \Pi^{od},\,
        (o, d) \in \mS.
\end{multline}

This set of constraints is of exponential size, but can be solved with a cutting-plane method that iterates between problems~\eqref{prob:maxan_benders_reformulation} and \eqref{prob:benders_sub_dual}. More specifically, we initialize Problem~\eqref{prob:maxan_benders_reformulation} without any of the constraints~\eqref{con:cut_plane}. We solve Problem~\eqref{prob:maxan_benders_reformulation} to obtain feasible $\bfx$ and $\bfz$ and trial values of $\tau^{od}$. For each $(o, d)\in \mS$, we then solve a Problem~\eqref{prob:benders_sub_dual} with the $\bfx$ and $\bfs$ and check if the trail value of $\tau^{od}$ and the optimal solution of Problem~\eqref{prob:benders_sub_dual} satisfy constraint~\eqref{con:cut_plane} or not. If not, the violated cut is added to Problem~\eqref{prob:maxan_benders_reformulation}. This process repeats until no new cut is added to Problem~\eqref{prob:maxan_benders_reformulation}. Although this Benders decomposition algorithm largely reduces the problem size and improves computational efficiency, it is insufficient to deal with our synthetic instances due to the widely acknowledged primal degeneracy issue for network flow problems \citep{magnanti1986tailoring} . We thus adopt a cut enhancement method and some acceleration strategies, which we describe in the next two sections, respectively.

\subsection{Pareto-Optimal Benders Cut}

We adapt the cut enhancement method proposed by \citet{magnanti1986tailoring} to generate pareto-optimal Benders cut by solving an auxiliary problem after a cut is identified by solving Problem~\eqref{prob:benders_sub_dual}. Let $(\bar{\bfx}, \bar{\bfz})$ denote a relative inner point of the feasible region specified by constraints \eqref{full:budget_edge}--\eqref{full:budget_node}, $\eta^{od}(\bfx, \bfz)$ denote the optimal value of Problem~\eqref{prob:benders_sub_dual} given $\bfx$ and $\bfz$. The Auxiliary problem associated with $(o, d) \in\mS$ is
\begin{subequations}
\begin{align}
    \underset{\boldsymbol{\lambda}, \boldsymbol{\theta}, \boldsymbol{\delta}, \boldsymbol{\pi} \geq 0}{\textrm{maximize}} \quad 
    & -\lambda_{d}^{od} + \lambda_{o}^{od}
        - \sum_{(i,j)\in \mE_h} 
            \bar{x}_{ij}\theta_{ij}^{od}
        - \sum_{i\in \mN_h}
            \sum_{(i,j)\in \mE^-_{h}(i)}
            \sum_{(w,l)\in \mE_{h}(i)}
            (\bar{x}_{wl} + \bar{z}_i) \delta_{ijwl}^{od} 
        - \sum_{(i, j)\in \mE} \pi^{od}_{ij}\\
    \textrm{subject to} \quad
    & -\lambda_{d}^{od} + \lambda_{o}^{od}
        - \sum_{(i,j)\in \mE_h} 
            x_{ij}\theta_{ij}^{od}
        - \sum_{i\in \mN_h}
            \sum_{(i,j)\in \mE^-_{h}(i)}
            \sum_{(w,l)\in \mE_{h}(i)}
            (x_{wl} + z_i) \delta_{ijwl}^{od} 
        - \sum_{(i, j)\in \mE} \pi^{od}_{ij}
        = \eta^{od}(\bfx, \bfz)
    \label{pareto:dual_optimal}\\
    & - \lambda_{j}^{od} 
        + \lambda_{i}^{od}
        - \mathds{1}\left[(i, j) \in \mE_h \right] \theta_{ij}^{od}
        - \mathds{1}(i\in \mN_h) 
            \sum_{(w,l)\in \mE_{h}(i)}
            \delta_{ijwl}^{od} 
        - \pi^{od}_{ij}
        \leq t_{ij},
        \quad \forall (i, j) \in \mE.
\end{align}
\label{prob:benders_pareto}
\end{subequations}
Problem \eqref{prob:benders_pareto} is different from problem \eqref{prob:benders_sub_dual} in that it has an additional constraint \eqref{pareto:dual_optimal} ensuring that the solution generated by problem \eqref{prob:benders_pareto} is optimal to problem \eqref{prob:benders_sub_dual}. We initialize the relative inner points as
\begin{subequations}
\begin{align*}
    & \bar{x}_{ij} = \min\{1,\  \frac{B_\text{edge}}{2|\mE_h|c_{ij}}\},
        \quad \forall (i, j)\in \mE_h, \\
    & \bar{z}_i = \min\{1,\ \frac{B_\text{node}}{2|\mN_h|b_i}\},
        \quad \forall i\in \mN_h.
\end{align*}
\label{eq:inner_point_init}
\end{subequations}
Through the solution process, every time an integral network design decision $(\bfx', \bfz')$ is found, we update the relative inner point as  
\begin{subequations}
\begin{align*}
    & \bar{x}_{ij} \gets \frac{1}{2} (\bar{x}_{ij} + x_{ij}'),
        \quad \forall (i, j)\in \mE_h, \\
    & \bar{z}_i \gets \frac{1}{2} (\bar{z}_{ij} + z_{ij}'),
        \quad \forall i\in \mN_h.
\end{align*}
\label{eq:inner_point_update}
\end{subequations}
This cut enhancement strategy requires longer time to generate a single cut as an auxiliary problem has to be solved. However, according to our computational experiments, it significantly reduces the number of cuts needed for solving the problem, and thus achieves shorter overall computation time compared to naive implementation of the Benders decomposition algorithm.

\subsection{Other Acceleration Strategies}

We adopt the following strategies to further accelerate the benders decomposition algorithm.
\begin{itemize}
    \item \textit{Initial Cut Generation.} Before solving problem \eqref{prob:maxan_benders_reformulation}, we apply the Benders decomposition algorithm to solve its linear-programming (LP) relaxation. Following \citet{bodur2017mixed}, we then add the cuts that are binding at the optimal solution of the LP relaxation to problem \eqref{prob:maxan_benders_reformulation}. These cuts help to obtain the LP-relaxation bound at the root node of the branch-and-bound tree.
    \item \textit{Flow Variable Reduction.} In Problem \eqref{prob:maxandp_piece_single}, routing variables are created for all $(o, d)\in \mS$ and all $(i, j) \in \mE$. However, given that a dummy low-stress link whose travel time is $\mT_2$ is added to connect each OD pair, edges that are far away from the origin and destination will not be used. Therefore, for each OD pair $(o, d) \in \mS$, we generate routing variables $x_{ij}^{od}$ only if the sum of travel time from $o$ to node $i$, travel time along edge $(i, j)$, and the travel time from node $j$ to $d$ is less than $T$. This pre-processing strategy significantly reduces the problem size.
    \item \textit{Design Variable Reduction}. In Problem \eqref{prob:maxan_benders_reformulation}, road design decisions for different edges are made separately. However, a more practical way of cycling infrastructure planning is to build continuous cycling infrastructure on road segments, each consisting of multiple edges. Such road segments can be identified through communication with transportation planners. We incorporate this consideration by replacing the road design variables $y_{ij}$ with $y_{p}$ where $p$ indicates the road segment that edge $(i, j)$ belongs to. In our computational and case studies, we group all edges between two adjacent intersections of arterial roads into one project, resulting in 84 and 1,296 projects in the synthetic and real networks, respectively. This preprocessing strategy reduces the number of binary decision variables. 
\end{itemize}

\subsection{Hyper-parameter tuning for the ML-augmented Model}

\subsection{Alternative Approach to Select \texorpdfstring{$k$}{Lg}}

In this section, we present a practical cross-validation approach to select the hyper-parameter $k$ for the $k$NN-augmented optimization model, which is summarized in Algorithm \ref{alg:follower_select_knn}.

\begin{algorithm}[!ht]
\small
\caption{A solution method using the $k$NN-augmented model}\label{alg:follower_select_knn}
\textbf{Input}: Width of the search window $\omega$; Number of leader decisions to sample $n_{d}$; Follower sample size $p$; Follower features $\{\bff^s\}_{s\in \mS}$. \\
\textbf{Output}: A leader solution $\hat{\bfx}^{k\text{NN}}$.
\begin{algorithmic}[1]
\State Randomly sample $n_{d}$ feasible leader solutions $\{\bfx_i \in \mX \}_{i=1}^{n_d}$.
\For{$i=1$ \textbf{to} $n_d$}
    \State Generate a dataset $\mD_i = \{\bff^s, G^s(\bfx_i)\}_{s\in \mS}$  .
    \State Perform a random train-test split to obtain $\mD_i^\text{train}$ and $\mD_i^\text{test}$ such that $|\mD_i^\text{train}| = p$.
    \For{$k=1$ \textbf{to} $p$}
        \State Build $k$NN model $P_{i, k}$ using $\mD^\text{train}_i$.
        \State Calculate out-of-sample loss $e_{i, k} = \frac{1}{|\mD_i^\text{test}|} \sum_{(\bff^s, G^s(\bfx_i))\in \mD_i^\text{test}} \left| P_{i, k}(\bff^s) - G^s(\bfx_i)\right|$. 
    \EndFor
\EndFor
\State Select the best $k^* \in \argmin_{k \in [p]} \left\{ \frac{1}{n_d} \sum^{n_d}_{i=1} e_{i, k}\right\}$.
\For{$k \in \mK := \left\{\max\{1,\, k^*-\omega\}, \dots, \min\{p,\, k^* + \omega\}\right\}$}
    \State Obtain leader's solution $\hat{\bfx}_k$ by solving problems~\eqref{prob:follower_select_knn} and \eqref{prob:non-parametric_approx} with $k$.
\EndFor
\State Select the best solution $\hat{\bfx}^{k\text{NN}} \in \argmin_{\bfx} \left\{F(\bfx) \,\middle|\, \bfx \in \{\hat{\bfx}_k\}_{k \in \mK} \right\}$.
\end{algorithmic}
\end{algorithm}

\section{Computational Study Details} \label{app:comp}

\subsection{The Synthetic Grid Network} \label{appsub:synthetic_network}
\begin{figure}[!ht]
    \centering
    \includegraphics[width=0.7\textwidth]{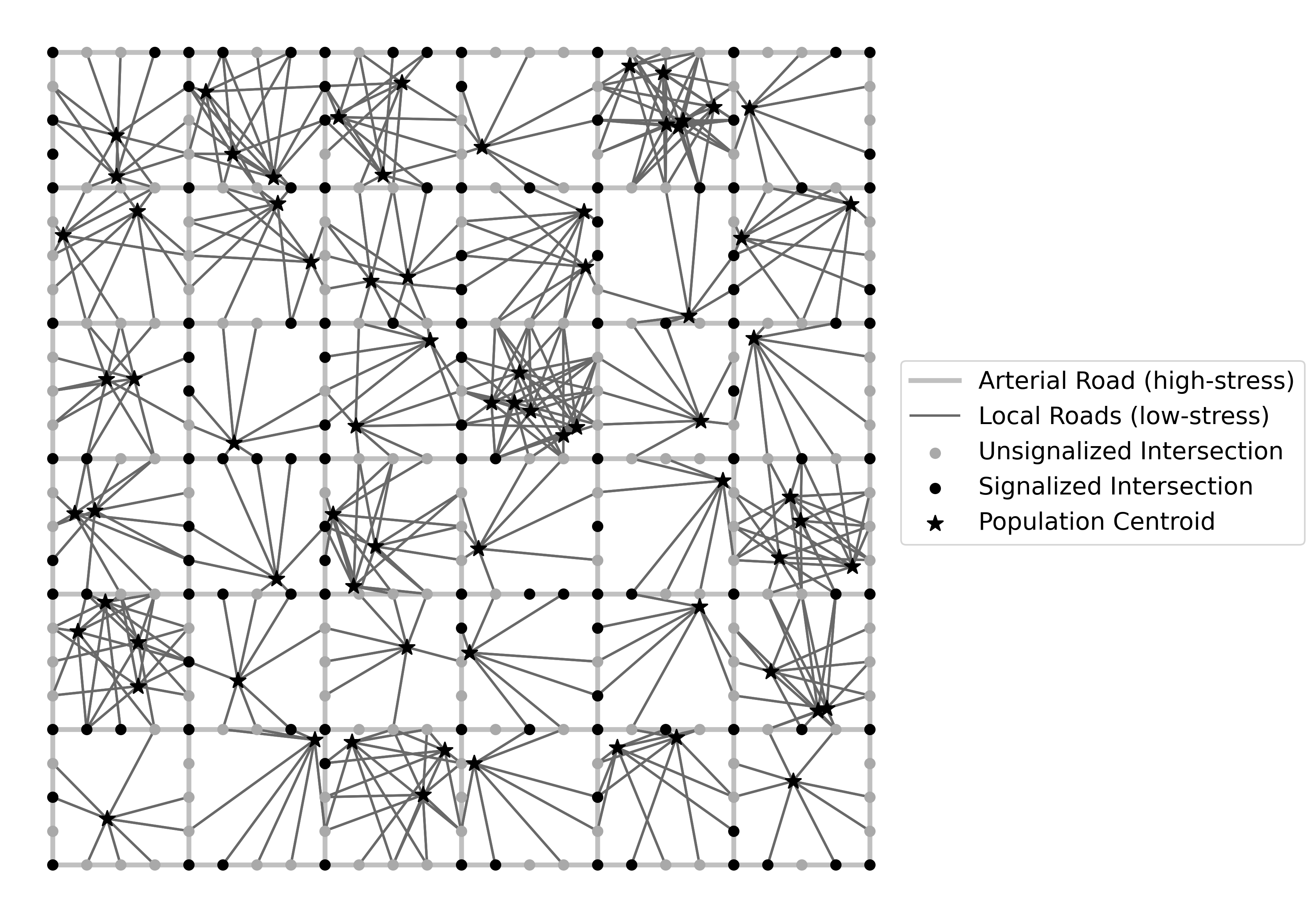}
    \caption{A synthetic grid network. Network components that are highlighted in black or dark grey constitute a low-stress network while others are high-stress.}
    \label{fig:synthetic_network}
\end{figure}

As presented in Figure \ref{fig:synthetic_network}, we create a synthetic network comprising a set of arterial roads, which are assumed to be high-stress, and local roads, which are assumed to be low-stress. The arterial roads constitute a 6x6 grid. Intersections of arterial roads are assumed to be signalized. On each arterial road segment, we place three nodes, each representing the intersection of the arterial road and a local road. Each of these intersections is assigned a traffic signal with a probability of 0.3. We generate 72 population centroids randomly distributed within the 36 major grid cells. Each centroid represents one origin and is a destination for all other centroids. We create low-stress edges that connect each centroid to 70\% of the intersections around the major grid cell in which that centroid is located. All edges are bidirectional. Each direction is assigned a travel time randomly distributed between 1 and 5. We use a constant travel speed of 1 to convert the generated travel time to distance. We consider all the OD pairs between which the shortest travel time on the overall network is less than 60, each assigned a weight uniformly distributed between 1 and 10. The network consists of 1,824 edges, 373 nodes, and 3,526 OD pairs. Arterial edges are grouped into 84 candidate projects. Setting the road design budget to 100, 300, and 500 roughly corresponds to selecting 5, 10, and 15 projects, respectively.


\subsection{Accessibility Calculations} \label{appsub:acc_measure}

\subsubsection{Location-based accessibility measures.} 

\begin{figure}[!ht]
    \centering
    \includegraphics[width=3in]{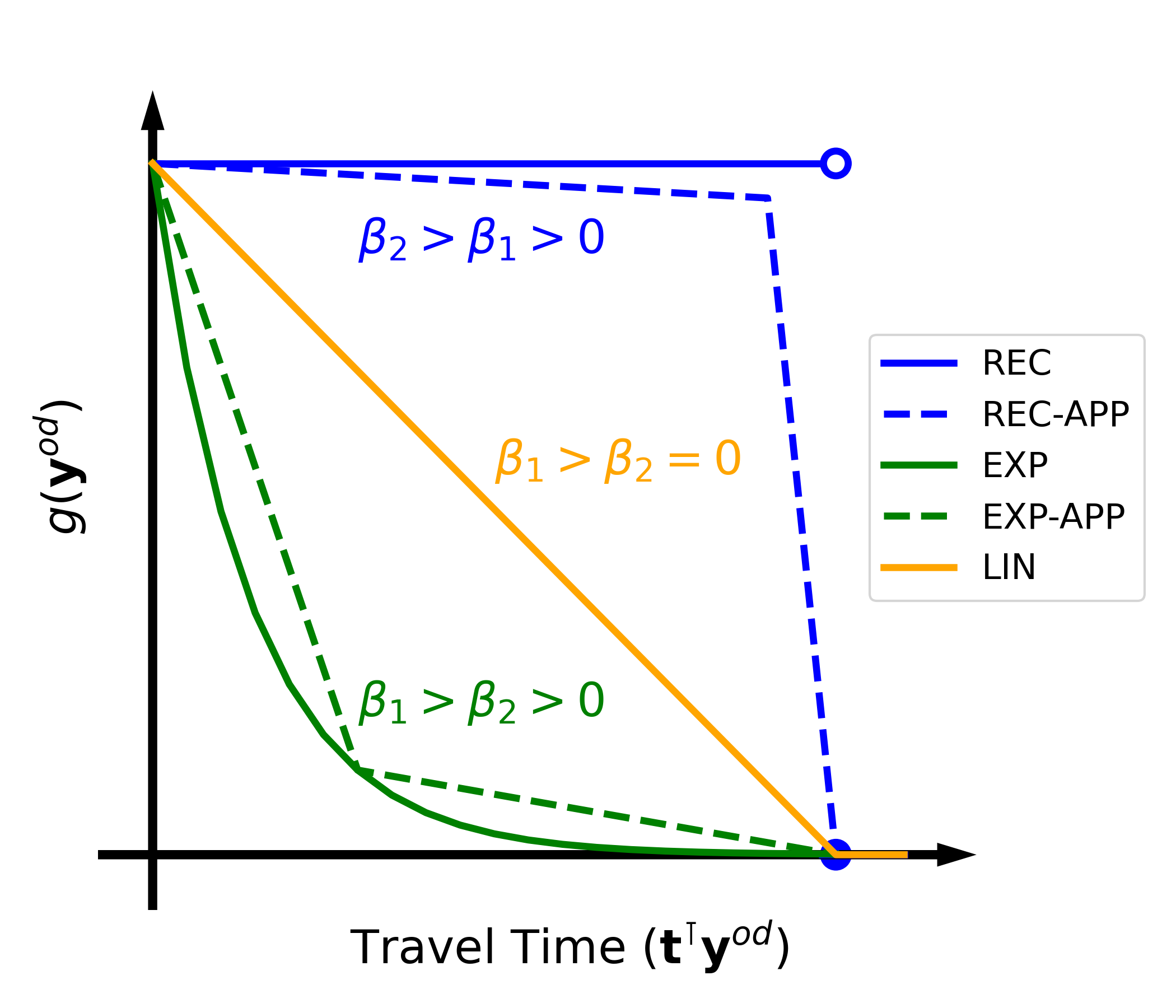}
    \caption{Negative exponential (EXP), rectangular (REC), and linear (LIN) impedance functions, and their piecewise linear approximations (APP).}
    \label{fig:impedance}
\end{figure}

We vary the parameters of the piecewise linear function $g$, mimicing three commonly used impedance function for location-based accessibility (Figure \ref{fig:impedance}):
\begin{enumerate}
    \item \textit{Negative exponential function}: $\beta_1=0.0375$, $\beta_2 = 0.00625$, $T_1=20$, $T_2=60$.
    \item \textit{Linear function}: $\beta_1=1/60$, $\beta_2 = 0$, $T_1=60$, $T_2=60$.
    \item \textit{Rectangular function}: $\beta_1=0.001$, $\beta_2 = 0.471$, $T_1=58$, $T_2=60$.
\end{enumerate}

\subsubsection{Utility-based accessibility measure.} We adopt the utility-based measure proposed by \citet{liu2021urban}. It requires as input a set of candidate routes for each OD pair and a constant $\alpha$ that reflects cyclists' preference between bike-lane continuity and bike-lane coverage along the routes. Following \citet{liu2021urban}, we set the value of $\alpha$ to $1.05$. We generate three candidate routes for each OD pair by solving three shortest path problems on the overall network using different edge travel cost. Specifically, we consider i) randomly generated travel time (Section \ref{fig:synthetic_network}), ii) Euclidean distance between the two ends, and iii) a uniform cost of 1 for all arterial roads and a uniform cost of 10 for other roads. The first and second definitions correspond to the goals of time minimization and distance minimization, both are commonly used by map software to generate route recommendations. The third definition reflects the preference for biking on major roads.

\subsection{Predictive Features} \label{appsub:feature}
Next, we present the implementation details of our representation learning framework (\ref{ec:rep_fea}), along with two baselines: TSP features (\ref{ec:tsp_fea}) and graph-theoretical features (\ref{ec:graph_fea}).

\subsubsection{Learned features.} \label{ec:rep_fea}
We apply the representation learning technique introduced in Section \ref{sec:repre_learning} to learn OD-pair features. Implementation details are illustrated as follows.

\begin{itemize}
    \item \textit{Relationship Graph Construction.} We sample $n_\text{sim}$ leader decisions for constructing the relationship graph. Ideally, the leader decisions should be sampled considering a specific design budget. However, the MaxANDP will be solved with various budgets in the computational studies (and in real-world transportation planning settings). Learning follower features multiple times may incur considerable computational burdens. Therefore, we consider learning features that are applicable to MaxANDP instances with different budgets by tweaking the leader decision sampling procedure. The high-level idea is to sample leader decisions with budgets randomly generated from a ``wide'' range that covers most budgets that may be used during the planning phase. As a result, the learned features contain information about follower similarity under various design budgets. Specifically, our leader-decision sampling procedure is parameterized by $\bar{P}$ and $\bar{Q}$ indicating, respectively, the maximum number of projects and the maximum number of nodes that can be selected in each sampled leader decision. Before generating a leader decision, we first randomly generate the number of projects and the number of nodes to be selected from intervals $[1, \bar{P}]$ and $[1, \bar{Q}]$, respectively. We then randomly select projects and nodes to form a leader decision. The procedure is summarized in Algorithm \ref{alg:sample_leader_dec}. In our computational studies, we set $\bar{P} = 25$, $\bar{Q} = 10$, and $n_\text{sim}$ = 5,000. We present computational results about the impact of $n_\text{sim}$ on feature quality in Section \ref{appsub:exp1_emb_params}.
    \item \textit{Follower Embedding.} In our computational study, we set $n_\text{walk} = 50$, $n_\text{length} = 20$, $\omega = 5$, and $\xi=16$. We investigate the impact of $\xi$ on feature quality in Section \ref{appsub:exp1_emb_params}.
\end{itemize}

\begin{algorithm}
\small
\caption{Sampling leader's decision for our cycling network design problem}\label{alg:sample_leader_dec}
\textbf{Input}: Number of decisions $n_{sim}$; Set of cycling infrastructure project $\mP$; Set of high-stress nodes $\mN_h$, Maximum number of infrastructure project selected $\bar{P}$; Maximum number of nodes selected $\bar{Q}$. \\
\textbf{Output}: A set of leader's decisions $\bar{\mX}$.
\begin{algorithmic}[1]
\State Initialize $\bar{\mX} = \{\}$.
\For{$i=1$ \textbf{to} $n_{sim}$}
    \State Generate $p \sim \textrm{Uniform}(1, \bar{P})$.
    \State Generate $q \sim \textrm{Uniform}(1, \bar{Q})$.
    \State Uniformly sample $\mP_i \subseteq \mP$ such that $|\mP_i| = p$.
    \State Uniformly sample $\mQ_i \subseteq \mN_h$ such that $|\mQ_i| = q$.
    \State Update $\bar{\mX} \gets \bar{\mX} \cup \{(\mP_i, \mQ_i)\}$
\EndFor
\end{algorithmic}
\end{algorithm}

\subsubsection{TSP features.} \label{ec:tsp_fea}
We adapt the features proposed by \citet{liu2021time} for predicting TSP objective values. Specifically, these features include
\begin{itemize}
    \item The coordinates of the origin.
    \item The coordinates of the destination.
    \item The Euclidean distance between the origin and the center of the grid.
    \item The Euclidean distance between the destination and the center of the grid.
    \item The Euclidean distance between the origin and the destination.
    \item The area of the smallest rectangular that covers both the origin and the destination.
    \item The travel time of the shortest path from the origin to the destination on the overall network.
\end{itemize}

\subsubsection{Graph theoretical features.} \label{ec:graph_fea}

\revision{We describe each follower $s\in\mS$ using the following node importance metrics in the relationship graph constructed in \ref{ec:rep_fea}.}

\begin{itemize}
    \item \revision{\textit{Degree centrality:} number of edges connected to node $s$.}
    \item \revision{\textit{Closeness centrality:} $1/\sum_{t\neq s} d(s, t)$, where $d(s, t)$ indicates the shortest path distance from node $s$ to node $t$.}
    \item \revision{\textit{Betweenness:} $\sum_{v \neq s \neq t} \sigma_{vt}(s)/\sigma_{vt}$, where $\sigma_{vt}$ is the total number of shortest paths from node $s$ to node $t$ and $\sigma_{vt}(s)$ is the number of those paths passing through node $s$.}
    \item \revision{\textit{Eigenvector centrality:} The $s$-th element of the left eigenvector associated with the eigenvalue of maximum modulus that is positive.}
    \item \revision{\textit{Current-flow centrality}: as defined by \citet{brandes2005centrality}.}
\end{itemize}

\subsection{Follower Sampling Methods} \label{appsub:follower_sampling}
We consider three follower sampling methods: i) vector assignment $p$-median sampling, ii) uniform sampling, and iii) $p$-center sampling. Uniform sampling selects each follower with a uniform probability, while the other two methods require solving an optimization problem. We next introduce the algorithms that we adopt to solve the optimization problems.

\subsubsection{Vector assignment \texorpdfstring{$p$}{Lg}-median sampling. }
We adapt the meta-heuristic proposed by \citet{boutilier2020ambulance} for solving the classical $p$-median problem to solve the vector assignment $p$-median problem introduced in Section \ref{subsec:practical}. Starting from a randomly generated initial solution $\mT$, the algorithm iterates between a ``swap'' phase and an ``alternation'' phase for $n_\text{iteration}$ iterations. In the swap phase, we create a solution in the neighborhood of $\mT$ by randomly removing one follower from $\mT$ and adding one follower from $\mS\backslash\mT$ to $\mT$. We update $\mT$ if the neighbor solution achieves a lower objective value for the vector assignment $p$-median problem. We perform at most $n_{swap}$ such swaps in each iteration. In the alternation phase, for each follower $t\in \mT$, we first find a follower set from $\mS$ whose $k$-nearest neighbor in $\mT$ includes $t$, and then solve an 1-median problem on this set. We solve in total $|\mT|$ 1-median problems and their optimal solutions form a new solution to the vector assignment $p$-median problem. We update $\mT$ if this solution achieves a lower objective value. The algorithm is summarized in Algorithm \ref{alg:meta_median}.

\begin{algorithm}[!ht]
\small
\caption{A Meta-Heuristic for Solving the Vector-Assignment $p$-Median Problem}\label{alg:meta_median}
\textbf{Input}: A set of followers $\mS$; Follower features $\{\bff^s\}_{s\in \mS}$; A distance metric in the feature space $d_{\mF}$; Number of followers to select $p$; Number of medians each node is assigned to $k$; Maximum number of iterations $n_\text{iteration}$; Number of swaps in each iteration $n_\text{swap}$; An oracle that calculates the objective value of the vector assignment $p$-median problem $O_{p-\text{median}}$ for any given solution; An oracle that finds the $k$-nearest neighbor of a follower in a given set $O_{k\text{NN}}$. \\
\textbf{Output}: A set of selected followers $\mT$.
\begin{algorithmic}[1]
\State Randomly sample $\mT\subseteq \mS$ such that $|\mT| = p$.
\State Initialize iteration counter $m_\text{interation} = 1$
\While{$m_\text{iteration} \leq n_\text{iteration}$}
    \For{$i=1$ \textbf{to} $n_\text{swap}$} \Comment{Swap}
        \State Randomly sample a follower $s \in \mS\backslash\mT$ and a follower $t\in \mT$.
        \State Create a follower set $\mT' = \mT\backslash \{t\} \cup \{s\}$
        \If{$O_{p-\text{median}}(\mT') < O_{p-\text{median}}(\mT)$} update $\mT\gets \mT'$ \EndIf
    \EndFor
\State Initialize a follower set $\mT'' = \{\}$ \Comment{Alternation}
\For{$t\in \mT$}
    \State Create a follower set $\mS^t = \left\{s \in \mS\backslash\mT \,\middle|\, t\in O_{k\text{NN}}(s, \mT) \right\}$.
    \State Solve an 1-median problem on $\mS^t$: $s'' \gets \argmin_{s\in \mS^t} \sum_{l\in \mS^t} d_\mF(\bff^s, \bff^l)$.
    \State Update $\mT'' \gets \mT''\cup \{s''\}$
\EndFor
\If{$O_{p-\text{median}}(\mT'') < O_{p-\text{median}}(\mT)$} update $\mT \gets \mT''$ \EndIf
\State Update $m_\text{iteration} \gets m_\text{iteration} + 1$
\EndWhile
\end{algorithmic}
\end{algorithm}

We note that the implementation of this meta-heuristic requires calculating the distance between every pair of followers in the feature space. For $\mS$ that is relatively small, we pre-calculate and store the distance matrix of the followers. However, storing the distance matrix (over 200 GB) in RAM is practically prohibitive for a large $\mS$ (e.g. the MaxANDP of Toronto's road network). To tackle the challenge, we calculate the distances on the fly when needed during the searching process. More specifically, we only calculate the distances from the current ``medians'' to other followers, resulting in a much smaller distance matrix that can be stored in RAM. A GPU with 24 GB of RAM (NVIDIA RTX6000) is employed to accelerate the distance-matrix calculation.

\subsubsection{\texorpdfstring{$p$}{Lg}-center sampling.}

We select followers by solving the follower problem
\begin{equation}
\label{prob:pcenter}
    \min_{\mT\subseteq\mS} 
    \left\{
        \max_{s\in \mS\backslash\mT} \min_{t\in \mT} d_\mF (\bff^s, \bff^t)
    \,\middle|\,
        |\mT| = p
    \right\}.
\end{equation}
We adapt a greedy heuristic to solve this problem. Specifically, we initialize the follower sample with a randomly selected follower $s \in \mS$. Next, we iteratively select one unselected follower and add it to the follower sample until $p$ followers have been selected. In each iteration, we first select calculate the shortest distance from each unselected follower to the follower sample, and then select the follower with the largest distance. This greedy heuristic generates 2-optimal solution for Problem~\eqref{prob:pcenter}. This is the best possible approximation that a heuristic algorithm can provide for the $p$-center problem because, for any $\delta < 2$, the existence of $\delta$-approximation implies $\mathcal{P}=\mathcal{NP}$ \citep{gonzalez1985clustering, hochbaum1985best}. To empirically improve the solution quality, we apply this algorithm for $n_\text{repeat}$ times and select the one that achieves the lowest objective value for Problem~\eqref{prob:pcenter}. For the computational experiments presented in Section \ref{sec:comp}, we set the value of $n_\text{repeat}$ to $200$. The heuristic is summarized in Algorithm \ref{alg:meta_center}.

\begin{algorithm}[!ht]
\small
\caption{A Greedy Heuristic for Solving the $p$-Center Problem}\label{alg:meta_center}
\textbf{Input}: A set of followers $\mS$; Follower features $\{\bff^s\}_{s\in \mS}$; A distance metric in the feature space $d_{\mF}$; Number of followers to select $p$; An oracle that calculates the objective value of the $p$-center problem $O_{p\text{-center}}$; Number of times the process is repeated $n_\text{repeat}$. \\ 
\textbf{Output}: A set of selected followers $\mT$.
\begin{algorithmic}[1]
\For{$i=1$ \textbf{to} $n_\text{repeat}$}
\State Randomly sample a follower $s \in \mS$.
\State Initialize follower sample $\mT_i = \{s\}$. 
\While{$|\mT_i| \leq p$}
    \State Calculate the distances to the selected set $d^s = \min_{t\in \mT_i} d_\mF (\bff^s, \bff^t)$ for all $s\in \mS\backslash\mT_i$.
    \State Select a follower $s' = \argmax_{s\in \mS\backslash\mT_i} d^s$.
    \State Update $\mT_i = \mT_i \cup \left\{s' \right\}$.
\EndWhile
\EndFor
\State Select $\mT = \argmin_{i=1 \in [n_\text{repeat}]} O_{p\text{-center}}(\mT_i)$
\end{algorithmic}
\end{algorithm}

\subsection{Euclidean norm-based Wasserstein scenario reduction.} \label{subapp:EW_sampling}
\revision{Following \citet{bertsimas2022optimization}, we employ a $k$-means-based heuristic to solve the $l_2$-norm Wasserstein scenario reduction problem. The steps are as follows:}
\begin{enumerate}
    \item \revision{Given followers $\{\bff^s\}_{s\in \mS}$, we apply the $k$-means algorithm to partition them into $k$ clusters.}
    \item \revision{For each cluster, we select the follower closest to the cluster center.}
\end{enumerate}

\subsection{DPCV scenario reduction.} \label{subapp:dpcv}
\revision{
This algorithm, adapted from \citep{dupavcova2003scenario}, adds one follower to our sample at the time until the number of followers selected is no less than our target. In each step, we select the follower that minimizes the $l_2$ Wasserstein distance between the sampled and unsampled followers. 
}

\subsection{Choosing \texorpdfstring{$\bar{L}$}{Lg}} \label{subapp:l_bar}

We propose a practical approach to iteratively search for an appropriate $\bar{L}$. The search starts from a given $L_0$, which is estimated using data associated with randomly generated leader decisions, and then gradually increases this value until the generated leader decision stops improving. The complete solution approach is presented as Algorithm \ref{alg:loss_reg}. 

\begin{algorithm}[!ht]
\small
\caption{A solution method using the parametric regression-augmented model}\label{alg:loss_reg}
\textbf{Input}: Step size $l_\text{step}$; Number of leader decisions to sample $n_{d}$; Follower sample size $p$; Follower features $\{\bff^s\}_{s\in \mS}$. \\
\textbf{Output}: A leader solution $\hat{\bfx}^\text{reg}$.
\begin{algorithmic}[1]
\State Randomly sample $n_{d}$ feasible leader solutions $\{\bfx_i \in \mX \}_{i=1}^{n_d}$.
\For{$i=1$ \textbf{to} $n_d$}
    \State Generate dataset $\mD_i = \{\bff^s, G^s(\bfx_i)\}_{s\in \mS}$.
    \State Randomly select a training set $\mD_i^\text{train}\subseteq \mD_i$ such that $|\mD_i^\text{train}| = p$.
    \State Train a prediction model $P_{i, k}\colon \bbR^\xi \rightarrow \bbR$ on $\mD^\text{train}_i$.
    \State Calculate the training loss $e_{i} = \frac{1}{|\mD_i^\text{train}|} \sum_{\bff^s, G^s(\bfx_i)\in \mD_i^\text{train}} \left| P_{i, k}(\bff^s) - G^s(\bfx_i)\right|$. 
\EndFor
\State Select a starting point $L_0 = \textrm{median} \left\{  e_1, e_2, \dots, e_{n_d}\right\}$.
\State Obtain $\mT$ by solving Problem \eqref{prob:follower_select_knn}.
\State Obtain an initial solution $\bfx_0$ by solving Problem \eqref{prob:reg_approx} with $L_0$ and $\mT$.
\State Initialize step counter $s = 1$
\Repeat
    \State Update $L_s = L_{s-1} + l_\text{step}$.
    \State Obtain $\bfx_s$ by solving Problem \eqref{prob:reg_approx} with $L_s$ and $\mT$.
\Until{$F(\bfx_s) > F(\bfx_{s-1})$}.
\State Select the best solution $\hat{\bfx}^\text{reg} = \bfx_{s-1}$.
\end{algorithmic}
\end{algorithm}

\subsection{Computational Setups}

All the algorithms were implemented using Python 3.8.3 on an Intel i7-8700k processor at 3.70 GHz and with 16GB of RAM. Optimization algorithms were implemented with Gurobi 9.1.2. The DeepWalk algorithm was implemented with Gensim 4.1.2. All the ML models for cycling accessibility prediction were implemented with Scikit Learn 1.0.2. 



\subsection{ML Model Implementation Details} \label{appsub:ML_param}
For the four ML models we consider in Section \ref{sec:comp}, we select the hyper-parameters (if any) based on the mean of median out-of-sample prediction performance over 1000 datasets. We note that we do not create a validation set because the goal is to achieve as good performance as possible on the out-of-sample follower set $\mS\backslash\mT$. The generalization of the ML models outside $\mS$ is not of interest in our study. The Linear regression does not involve any hyper-parameter. The neighborhood sizes of $kNN$, the regularization factors of the lasso regression and ridge regression are summarized in Tables \ref{tab:knn_params}--\ref{tab:ridge_params}, respectively.

\begin{table}[!ht]
\centering
\caption{Neighborhood sizes of $k$-nearest neighbor regression.} \label{tab:knn_params}
\scriptsize
\begin{tabular}{@{}lrrrrr@{}}
\toprule
                          & \multicolumn{1}{l}{}       & \multicolumn{4}{l}{Accessibility Measure}                                                            \\ \cmidrule(l){3-6} 
Feature                   & \multicolumn{1}{l}{Budget} & \multicolumn{1}{l}{EXP} & \multicolumn{1}{l}{LIN} & \multicolumn{1}{l}{REC} & \multicolumn{1}{l}{UT} \\ \midrule
\multirow{3}{*}{Learning} & 100                        & 1                       & 1                       & 1                       & 1                      \\
                          & 300                        & 1                       & 1                       & 1                       & 1                      \\
                          & 500                        & 1                       & 1                       & 1                       & 1                      \\ \midrule
\multirow{3}{*}{TSP}      & 100                        & 1                       & 1                       & 1                       & 1                      \\
                          & 300                        & 1                       & 1                       & 1                       & 1                      \\
                          & 500                        & 1                       & 1                       & 1                       & 1                      \\ \bottomrule
\end{tabular}
\end{table}

\begin{table}[!ht]
\centering
\caption{Regularization parameters of lasso regressions.} \label{tab:lasso_params}
\scriptsize
\begin{tabular}{@{}lrrrrr@{}}
\toprule
                          & \multicolumn{1}{l}{}       & \multicolumn{4}{l}{Accessibility Measure}                                                            \\ \cmidrule(l){3-6} 
Feature                   & \multicolumn{1}{l}{Budget} & \multicolumn{1}{l}{EXP} & \multicolumn{1}{l}{LIN} & \multicolumn{1}{l}{REC} & \multicolumn{1}{l}{UT} \\ \midrule
\multirow{3}{*}{Learning} & 100                        & 0.004                   & 0.008                    & 0.020                  & 0.020                   \\
                          & 300                        & 0.004                   & 0.008                    & 0.010                  & 0.020                   \\
                          & 500                        & 0.004                   & 0.006                    & 0.010                  & 0.020                   \\ \midrule
\multirow{3}{*}{TSP}      & 100                        & 0.006                   & 0.010                    & 0.040                  & 0.020                   \\
                          & 300                        & 0.008                   & 0.010                    & 0.030                  & 0.030                  \\
                          & 500                        & 0.007                   & 0.010                    & 0.020                  & 0.040                  \\ \bottomrule
\end{tabular}
\end{table}

\begin{table}[!ht]
\centering
\scriptsize
\caption{Regularization parameters of ridge regressions.} \label{tab:ridge_params}
\begin{tabular}{@{}lrrrrr@{}}
\toprule
                          & \multicolumn{1}{l}{}       & \multicolumn{4}{l}{Accessibility Measure}                                                            \\ \cmidrule(l){3-6} 
Feature                   & \multicolumn{1}{l}{Budget} & \multicolumn{1}{l}{EXP} & \multicolumn{1}{l}{LIN} & \multicolumn{1}{l}{REC} & \multicolumn{1}{l}{UT} \\ \midrule
\multirow{3}{*}{Learning} & 100                        & 50                      & 60                      & 40                      & 30                    \\
                          & 300                        & 50                      & 40                      & 10                      & 60                    \\
                          & 500                        & 30                      & 20                      & 10                      & 100                   \\ \midrule
\multirow{3}{*}{TSP}      & 100                        & 260                     & 240                     & 140                     & 150                   \\
                          & 300                        & 210                     & 170                     & 170                     & 150                    \\
                          & 500                        & 150                     & 150                     & 150                     & 100                    \\ \bottomrule
\end{tabular}
\end{table}

\subsection{Additional Results on the Impact of Hyperparameters on Feature Quality} \label{appsub:exp1_emb_params}
We focus on two hyperparameters that may affect the quality of the REP features: i) feature dimensionality $\xi$, and ii) the number of leader decisions to sample $n_\text{sim}$. The idea is to first use a relatively large $n_\text{sim}$, which makes sure the embedding algorithm is well-informed about the relationship between followers, to find the smallest $\xi$ that supports ML models to achieve ``good'' prediction performance. We want $\xi$ to be small because it helps to reduce the size of the ML-augmented model. Once a $\xi$ is chosen, we then search for a small $n_\text{sim}$ that makes the learned features perform well. We want $n_\text{sim}$ to be small because it reduces the computational efforts in constructing the follower relationship graph. We focus on the prediction performance of $k$NN using UNI samples because $k$NN constantly achieves the best prediction performance among all ML models considered.

\subsubsection{The impact of \texorpdfstring{$\xi$}{Lg}.} We follow the convention to set the value of $\xi$ to the powers of two. We vary $\xi$ in $\{2, 4, 8, 16\}$ because the synthetic network has 3,526 followers and the smallest training sample considered is 1\% of them, corresponding to 35 followers. Training linear regression models using 35 data points and features of 32 dimensions or more may lead to serious overfitting issues. We set $n_\text{sim}$ to 5,000 when choosing $\xi$.

\begin{figure}[!ht]
    \centering
    \includegraphics[width=.8\textwidth]{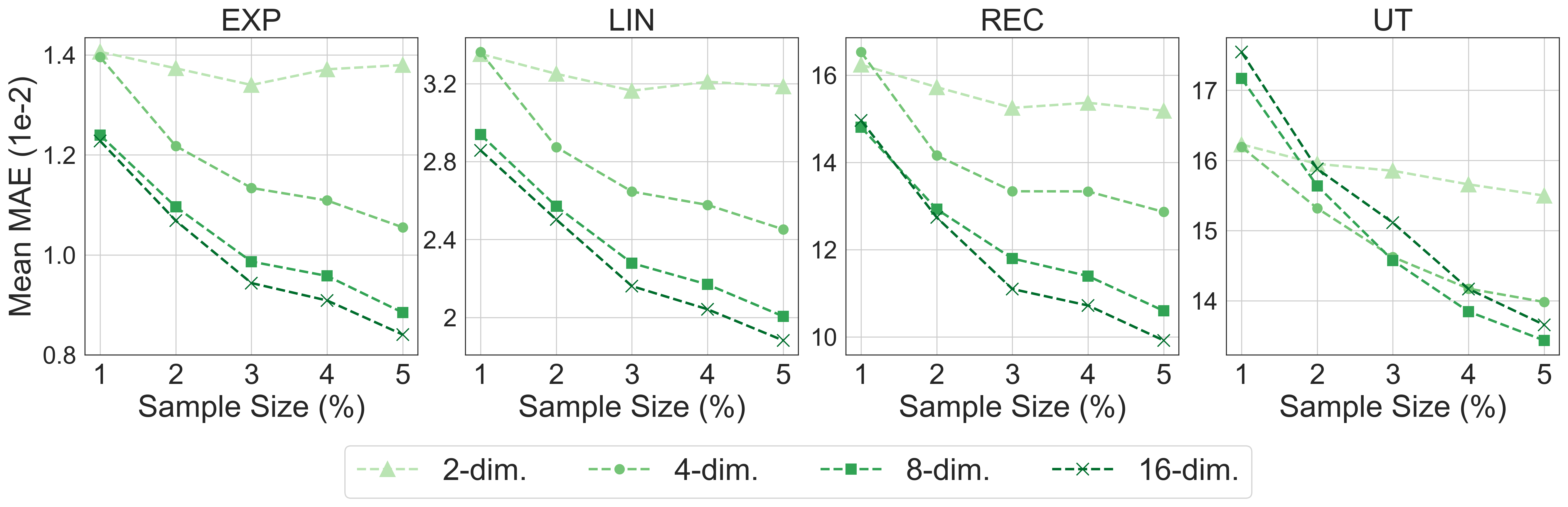}
    \caption{Out-of-sample prediction performance of $k$NN using REP features of 2, 4, 8, and 16 dimensions.}
    \label{fig:impact_xi}
\end{figure}

Figure \ref{fig:impact_xi} summarizes the predictivity of REP features of different dimensions. For location-based accessibility measures, increasing the number of dimensions leads to lower out-of-sample prediction error. For utility-based accessibility measures, REP features of four dimensions achieve the lowest error when the sample size is 1\% or 2\% of the original follower set, while 8- and 16-dimensional REPs become more predictive with larger training samples. For convenience, we choose to set $\xi = 16$ for all accessibility measures. However, based on the results presented in Figure \ref{fig:impact_xi}, carefully tuning $\xi$ for different problems may improve ML models' out-of-sample prediction accuracy.

\subsubsection{The impact of \texorpdfstring{$n_\text{sim}$}{Lg}.} 
\begin{figure}[!ht]
    \centering
    \includegraphics[width=.8\textwidth]{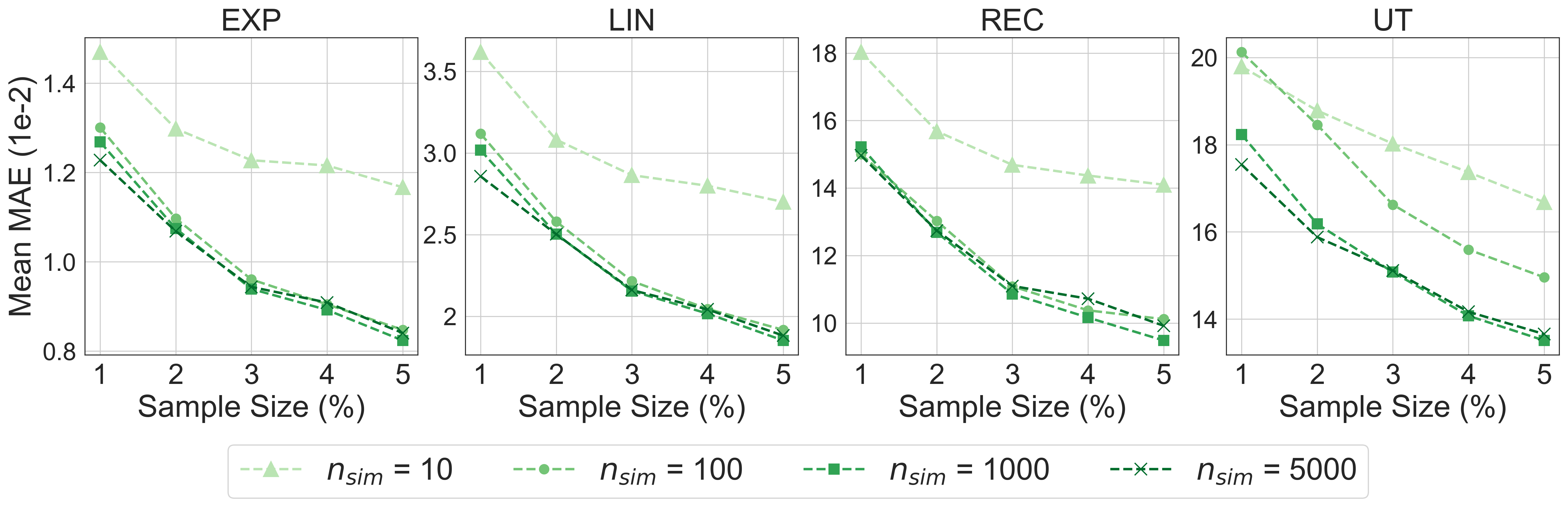}
    \caption{Out-of-sample prediction performance of $k$NN using REP features learned using 10, 100, 1000, and 5000 sampled leader decisions.}
    \label{fig:impact_nsim}
\end{figure}

We then fix $\xi=16$ and vary $n_\text{sim}$ in \{10, 100, 1,000, 5,000\}. Figure \ref{fig:impact_nsim} summarizes the predictivity of REP features learned with different $n_{sim}$. In general, considering more leader decisions lead to better prediction performance of the REP features, but the marginal improvement decreases as $n_\text{sim}$ increases. For location-based accessibility measures, the improvement from $n_\text{sim}=100$ to $n_\text{sim}=$5,000 is negligible compared to the improvement from $n_\text{sim}=10$ to $n_\text{sim}=100$. For the utility-based measure, REP features learned with $n_\text{sim}$=1,000 is as predictive as REP features learned with $n_\text{sim}=$ 5,000.

\section{Case Study Details}  \label{app:case}

\subsection{Network Construction and Pre-Processing} \label{appsub:preprocessing}
We retrieve Toronto centerline network from Toronto Open Data Portal \citep{trt2020open}. We follow the following steps to process the network data.
\begin{enumerate}
    \item We remove roads where cycling is physically impossible or legally prohibited, including ``highway'', ``highway ramp'',  ``railway'',  ``river'', ``hydro line'', ``major shoreline'', ``major shoreline (land locked)'', ``geostatistical line'', ``creek/tributary'', ``ferry lane'' per the City's definition. 
    \item We remove all redundant nodes. A node is considered redundant if it is not an end of a road or an intersection of three or more edges. These nodes are included in the original network to depict the road shape, but are unnecessary from a modeling perspective. 
    \item We replace local roads with low-stress edges that connect DA centroids to intersections along arterial roads. We solve a shortest-path problem from each DA centroid to each intersection located on its surrounding arterial roads using the low-stress network. If a low-stress path is found, we add a bi-directional low-stress edge that connects the DA centroid and the intersection and set its travel time to the travel time along the path. All local roads are then removed because we do not consider building new cycling infrastructure on local roads, and the role of local roads in our problem is to connect DA centroids to arterial roads, which can be served by the added artificial edges. The node and edge removal procedures are illustrated in Figure \ref{fig:remove_redundant_components}.\begin{figure}[!ht]
    \centering
    \includegraphics[width=.8\textwidth]{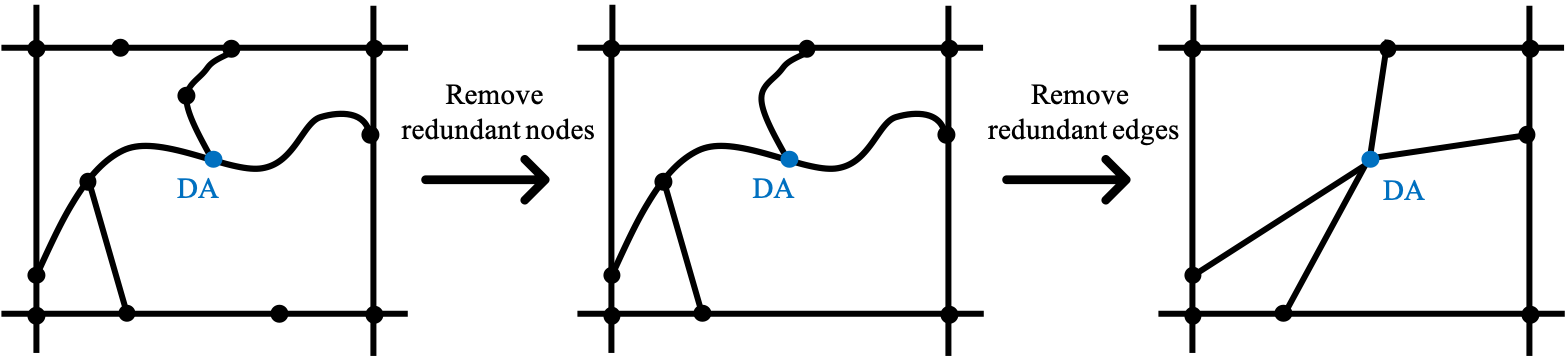}
    \caption{The procedures of removing redundant nodes and edges.}
    \label{fig:remove_redundant_components}
    \end{figure}
    \item We group arterial edges to form candidate projects. A candidate project is defined as a continuous road segment that connects two adjacent intersections of arterial roads. Such a road segment may be represented as multiple arterial edges in the network due to the presence of arterial-local intersections. Grouping these edges together allows us to create fewer road design variables. The average length of the projects is 1.48 km.
    \item Each DA is represented by its geometric centroid and is manually connected to its nearest point in the cycling network using a low-stress edge with zero travel time.
    \item We discard all the OD pairs that are currently connected via a low-stress path because building new cycling infrastructure will not affect the accessibility of these OD pairs.
\end{enumerate}

\subsection{Follower Embedding Details}
We apply similar procedures as introduced in Section \ref{appsub:feature} to learn follower features. We highlight the key difference in each step as follows.
\begin{itemize}
    \item \textit{Relationship Graph Construction}. We sample $n_\text{sim}=$10,000 leader decisions with $\bar{P} = 300$ and $\bar{Q} = 100$. We note that 149,496 (11.3\%) of the 1,327,849 OD pairs have zero accessibility under all sampled leader decisions. As a result, their similarities to other OD pairs are all zero according to the adopted similarity measure. These OD pairs are mostly outside downtown Toronto, where the road networks are highly stressful. Connecting these OD pairs with low-stress routes requires constructing a large amount of new cycling infrastructure, which is beyond the considered budget (100 km). We choose to exclude these OD pairs from OD pair embedding and thus exclude them from OD pair selection and ML-augmented model to avoid the computational burdens of sampling more leader decisions. However, we do take these OD pairs into consideration when evaluating leader decisions.
    \item \textit{Follower Embedding}. We set $n_\text{walk} = 50$, $n_\text{length} = 50$, $\omega = 5$, and $\xi=32$.
\end{itemize}

\subsection{Computational Setups} \label{appsub:trt_comp}
Optimization algorithms were implemented with Python 3.8.3 using Gurobi 9.1.2 on an Intel i7-8700k processor at 3.70 GHz and with 16GB of RAM. The heuristic for solving vector-assignment p-median problem is accelerated with h an NVIDIA P100 GPU. The DeepWalk algorithm was implemented with Gensim 4.1.2. 

For optimal network expansions, we use follower samples of 2,000 followers (OD pairs), and solve the $k$NN-augmented model with them. We use the $k$NN-augmented model because the network design budget ($\leq$ 100 km) falls into the small budget regime, where the $k$NN-augmented model generally outperforms the linear regression-augmented as presented in Section \ref{subsec:exp_opt}. We set the solution time limit to 3 hours for all optimization models. The greedy network expansion is implemented using the same machine as used by the optimization models and is parallelized using eight threads.

\subsection{Comparison between Greedy and Optimal Expansions} \label{appsub:greedy_opt_comp}

\begin{figure}[!th]
    \centering
    \includegraphics[width=0.6\textwidth]{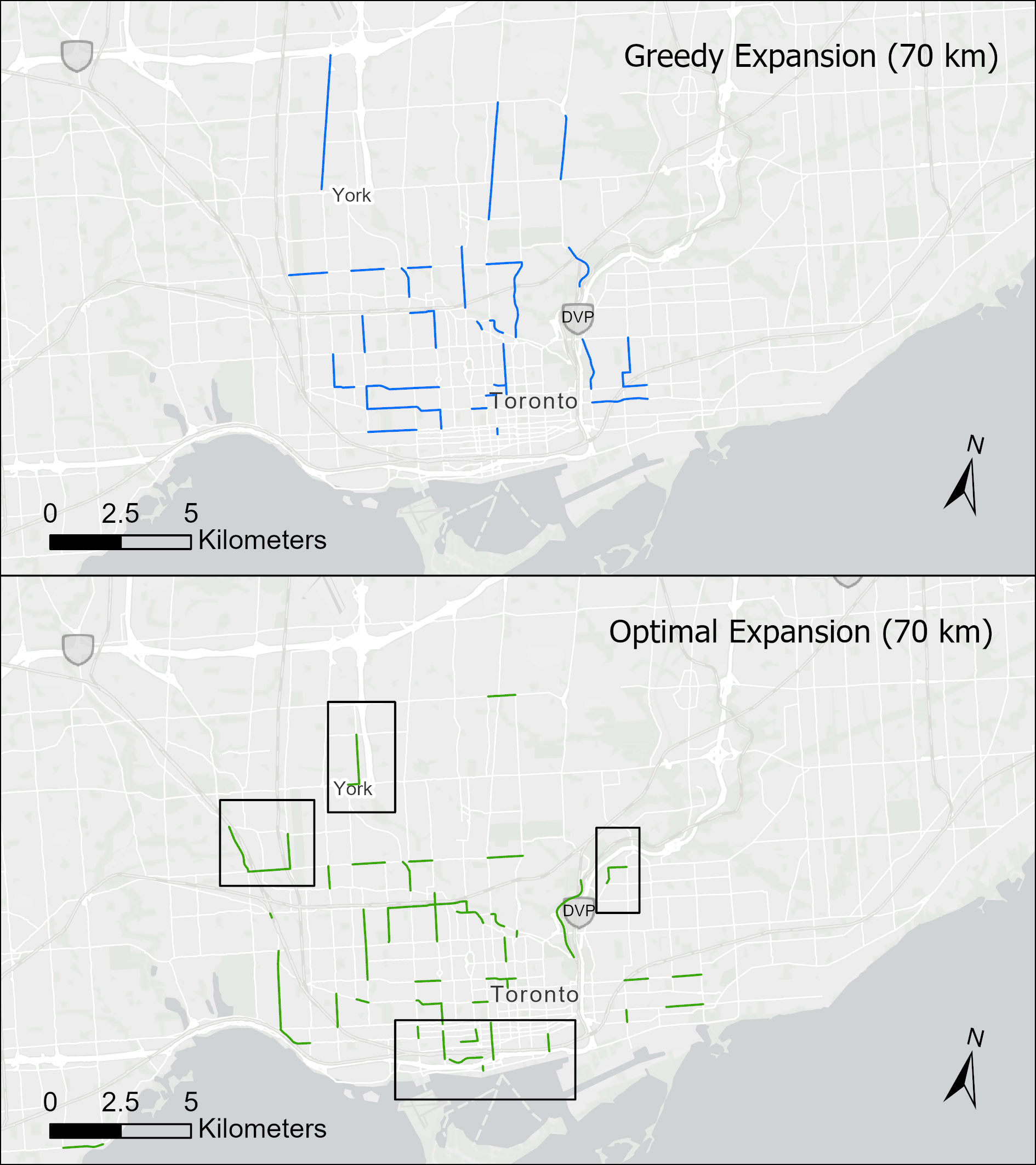}
    \caption{Greedy and optimal expansions given a road design budget of 70 km.}
    \label{fig:greedy_opt_compare}
\end{figure}

Figure \ref{fig:greedy_opt_compare} presents the cycling infrastructure projects selected by the greedy heuristic and our approach given a road design budget of 70 km, as an example. The optimal expansion is 11.2\% better than the heuristic expansion as measured by the improvement in Toronto's total low-stress cycling accessibility. Both algorithms choose many cycling infrastructure projects in the downtown core area, where a well-connected low-stress cycling network has already been constructed, and where job opportunities are densely distributed. These projects connect many DAs to the existing cycling network and thus grant them access to job opportunities via the existing network. However, unlike the greedy heuristic that spends almost all the road design budgets to expand the existing network, our approach identifies four groups of projects that are not directly connected to the existing network (as highlighted by the black frames in Figure \ref{fig:greedy_opt_compare}). The greedy heuristic does not select these projects because they have little impact on the total cycling accessibility if constructed alone. However, when combined, these projects significantly improve the accessibility of their surrounding DAs by breaking the high-stress barriers between low-stress cycling islands. 

\end{document}